\documentclass[a4paper,10pt]{article}
\usepackage{stmaryrd}
\usepackage{amsfonts}
\usepackage{bbm}
\usepackage{amscd}
\usepackage{mathrsfs}
\usepackage{latexsym,amssymb,amsmath,amscd,amscd,amsthm,amsxtra}
\usepackage[dvips]{graphicx}
\usepackage[utf8]{inputenc}
\usepackage[T1]{fontenc}
\usepackage{lmodern}
\usepackage{amssymb}
\usepackage[all]{xy}
\usepackage{nicefrac,mathtools,enumitem}
\usepackage[numbers,sort&compress]{natbib}
\usepackage{microtype}
\usepackage{xcolor}

\newtheorem{thm}{Theorem}[section]
\newtheorem{lem}[thm]{Lemma}
\newtheorem{cor}[thm]{Corollary}
\newtheorem{pro}[thm]{Proposition}
\newtheorem{ex}[thm]{Example}
\newtheorem{rmk}[thm]{Remark}
\newtheorem{defi}[thm]{Definition}

\newcommand {\emptycomment}[1]{}

\newcommand{\lon }{\,\rightarrow\,}
\newcommand{\be }{\begin{equation}}
\newcommand{\ee }{\end{equation}}

\newcommand{\pf}{\noindent{\bf Proof.}\ }

\newcommand{\g}{\frkg}

\newcommand{\AD}{\mathfrak{ad}}

\newcommand{\aaa}{\phi_{A}}
\newcommand{\aee}{\phi_{E}}

\newcommand{\aab}{A}
\newcommand{\aar}{\Gamma(A)}

\newcommand{\aap}{\varphi^*}
\newcommand{\aad}{\phi_{A}^{\dag}}
\newcommand{\aaq}{\Gamma(A^*)}
\newcommand{\aadd}{(\phi_{A}^{\dag})^{-1}}
\newcommand{\aae}{E}
\newcommand{\aaf}{\Gamma(\wedge^kA)}
\newcommand{\aag}{\Gamma(\wedge^kA^*)}
\newcommand{\aaee}{\Gamma(E)}
\newcommand{\ddd}{\mathcal{D}}

\newcommand{\Real}{\mathbb R}

\newcommand{\huaB}{\mathcal{B}}

\newcommand{\huaA}{\mathcal{A}}
\newcommand{\huaL}{\mathcal{L}}

\newcommand{\huaF}{\mathcal{F}}

\newcommand{\huaV}{\mathcal{V}}

\newcommand{\huaD}{\mathcal{D}}

\newcommand{\CWM}{C^{\infty}(M)}

\newcommand{\frkd}{\mathfrak d}

\newcommand{\frkg}{\mathfrak g}

\newcommand{\frkA}{\mathfrak A}

\def\qed{\hfill ~\vrule height6pt width6pt depth0pt}

\newcommand{\half}{\frac{1}{2}}

\newcommand{\br}[1]{   [ \cdot,    \cdot  ]_\frkg   }

\newcommand{\dev}{\mathfrak{D}}

\newcommand{\dM}{\mathrm{d}}

\newcommand{\Der}{\mathrm{Der}}

\newcommand{\Ad}{\mathrm{Ad}}

\newcommand{\gl}{\mathfrak {gl}}

\newcommand{\idd}{\mathrm{id}}

\newcommand{\ad}{\mathrm{ad}}

\begin{document}
\title{
{Hom-Lie algebroids, Hom-Lie bialgebroids and Hom-Courant algebroids
\thanks
 {
Research supported by NSFC (11471139) and NSF of Jilin Province (20140520054JH).
 }
} }
\author{Liqiang Cai$^1$, Jiefeng Liu$^2$ and Yunhe Sheng$^3$  \\
$^1$School of Mathematics and Statistics, Henan University,\\
Kaifeng 475000, Henan, China\\
$^2$Department of Mathematics, Xinyang Normal University,\\ Xinyang 464000, Henan, China\\
$^3$Department of Mathematics, Jilin University,\\
 Changchun 130012, Jilin, China
\\\vspace{3mm}
Email: cailq13@mails.jlu.edu.cn,~jfliu13@mails.jlu.edu.cn,~shengyh@jlu.edu.cn }

\date{}


\maketitle
\begin{abstract}
In this paper, first we modify the definition of a Hom-Lie algebroid introduced by Laurent-Gengoux and Teles and give its equivalent dual description. Many results that parallel to  Lie algebroids  are given. In particular, we give the notion of a Hom-Poisson manifold and show that there is a Hom-Lie algebroid structure on the pullback of the cotangent bundle of a Hom-Poisson manifold. Then we give the notion of a Hom-Lie bialgebroid, which is a natural generalization of a purely Hom-Lie bialgebra and a Lie bialgebroid. We show that the base manifold of a Hom-Lie bialgebroid  is a Hom-Poisson manifold. Finally, we introduce the notion of a Hom-Courant algebroid and show that  the double of a Hom-Lie bialgebroid is a Hom-Courant algebroid. The underlying algebraic structure of a Hom-Courant algebroid is a Hom-Leibniz algebra, or a Hom-Lie 2-algebra.
\end{abstract}

\tableofcontents

\section{Introduction}
In the study of $\sigma$-derivations of an associative algebra, Hartwig, Larsson
and Silvestrov introduced the notion of a Hom-Lie algebra in \cite{HLS}. Some $q$-deformations of the Witt and the Virasoro algebras have the
structure of a Hom-Lie algebra (\cite{HLS,Hu, Ka, LD1,LD2}).
Then in \cite{MS2}, Makhlouf and Silvestrov modified the definition of a Hom-Lie algebra. In Makhlouf and Silvestrov's new definition, a Hom-Lie algebra $(\g,[\cdot,\cdot]_\g,\phi_\g)$ is a nonassociative algebra $(\g,[\cdot,\cdot]_\g)$ together with an algebra homomorphism $\phi_\g:\g\longrightarrow \g$ such that $[\cdot,\cdot]_\g$ is skew-symmetric and the following Hom-Jacobi identity holds:
$$
  [\phi_\g(x),[y,z]_\g]_\g+[\phi_\g(y),[z,x]_\g]_\g+[\phi_\g(z),[x,y]_\g]_\g=0,\quad \forall x,y,z\in \g.
$$
On the set of $(\sigma,\sigma)$-derivations of a commutative algebra, there is a natural Hom-Lie algebra structure in the sense of Makhlouf
and Silvestrov (\cite{CS-bialgebra}). In the sequel, we only consider  Makhlouf
and Silvestrov's Hom-Lie algebra. Recently, the interest in Hom-algebras grows again due to the work of geometrization of Hom-Lie algebras (\cite{LGT}), quantization of Hom-Poisson structures (\cite{BEM-quantization}) and integration of Hom-Lie algebras (\cite{LGMT}).

The notion of a Lie algebroid was introduced
by Pradines in 1967, which is a generalization of Lie algebras and
tangent bundles. See \cite{General theory of Lie groupoid and Lie
algebroid} for general theory about Lie algebroids. They play
important roles in various parts of mathematics. In \cite{LGT},  Laurent-Gengoux and Teles introduced the notion of a Hom-Lie algebroid with the help of Hom-Gerstenhaber algebra. They showed that there is a one-to-one correspondence between Hom-Gerstenhaber algebra structures on $\Gamma(\wedge^\bullet \mathscr{A})$ and Hom-Lie algebroid structures on a vector bundle $\mathscr{A}$. The work in \cite{LGT} could serve to broaden the scope in important directions as the future
mathematical landscape continues to develop.

The first aim of this paper is to develop theories that parallel to Lie algebroids for Hom-Lie algebroids. First we introduce the notion of a Hom-bundle to make the definition of a Hom-Lie algebroid concise. There is a slight difference between our Hom-Lie algebroid and the one introduced in \cite{LGT}. We give the dual description of a Hom-Lie algebroid using twisted differential graded commutative algebra. On the pullback bundle of a Lie algebroid, there is a Hom-Lie algebroid structure. In particular, for any diffeomorphism $\varphi:M\longrightarrow M$, the pull back $\varphi^!TM$ of the tangent Lie algebroid $TM$ is naturally a Hom-Lie algebroid. On one hand, this can be viewed as the fundamental example of Hom-Lie algebroids and verifies  the validity of our definition of a Hom-Lie algebroid. On the other hand, this shows the naturality of the Hom-Lie algebroid structure. We show that the anchor map of a Hom-Lie algebroid is a Hom-Lie algebroid homomorphism to  $\varphi^!TM$. 
Using the Hom-Lie algebroid structure on $\varphi^!TM$, we introduce the notion of a Hom-Poisson tensor, which is equivalent to a purely Hom-Poisson algebra structure on $\CWM$. This is nontrivial. Without the discovery of the Hom-Lie algebroid $\varphi^!TM$, the definition of a Hom-Poisson manifold could not be given so elegantly. Finally, we show that  there is  a Hom-Lie algebroid structure on $\varphi^!T^*M$ associated to a Hom-Poisson manifold.

The second aim of this paper is to study the   bialgebroid theory for Hom-Lie algebroids. We give the definition of a Hom-Lie bialgebroid using the Hom-Gerstenhaber algebra structure given in \cite{LGT} and the differential operator given in this paper.  We show that a Hom-Poisson manifold gives rise to a Hom-Lie bialgebroid naturally and the base manifold of a Hom-Lie bialgebroid is a Hom-Poisson manifold. This generalizes the classical results about Poisson manifolds and Lie bialgebroids. Then we introduce the notion of a Hom-Courant algebroid and show that on the double of a Hom-Lie bialgebroid, there is naturally a Hom-Courant algebroid structure. Finally, we investigate the algebraic structure underlying a Hom-Courant algebroid and show that a Hom-Courant algebroid gives rise to a Hom-Lie 2-algebra, which is the categorification of a Hom-Lie algebra (\cite{SC-2-algebra}).

Note that a Hom-Lie algebroid is not only a formal generalization of a Lie algebroid. It has its own geometric meaning and fruitful examples, such as aforementioned the Hom-Lie algebroid $\varphi^!TM$ and the Hom-Lie algebroid associated to a Hom-Poisson manifold. This makes our result full of interest.

The paper is organized as follows. In Section 2, we give a review of representations of Hom-Lie algebras and purely Hom-Lie bialgebras. In Section 3, first we give the notion of a Hom-bundle and provide several examples. Then in Subsection 3.1, we give the definition of a Hom-Lie algebroid and give its dual description. In Subsection 3.2, we give the Hom-Lie algebroid $\varphi^!TM$ and show that the anchor of a Hom-Lie algebroid is a Hom-Lie algebroid homomorphism to $\varphi^!TM$. In Subsection 3.3, we give the formula of a Lie derivative for a Hom-Lie algebroid and provide many useful formulas. In Subsection 3.4, we introduce  the notion of a  Hom-Poisson manifold and show that there is a natural Hom-Lie algebroid structure on $\varphi^!T^*M$ associated to a Hom-Poisson manifold.   In Section 4, we introduce the notion of a Hom-Lie bialgebroid and show that the base manifold of a Hom-Lie bialgebroid is a Hom-Poisson manifold. In Section 5, we introduce the notion of a Hom-Courant algebroid, and show that the double of a Hom-Lie bialgebroid is a Hom-Courant algebroid and a Hom-Courant algebroid gives rise to a Hom-Lie 2-algebra.

\section{Preliminaries}


Let $\frkA$ be an associative algebra,   $\sigma$ and $\tau$  two   algebra endomorphisms on $\frkA$. A {\bf $(\sigma,\tau)$-derivation} on $\frkA$ is a linear map $D:\frkA\longrightarrow \frkA$ such that
$$
  D(ab)=D(a)\tau(b)+\sigma(a)D(b),\quad \forall~a,b\in\frkA.
$$
 The set of all $(\sigma,\tau)$-derivations on $\frkA$ is denoted by $\Der_{\sigma,\tau}(\frkA).$

 The notion of a $(\sigma,\tau)$-differential graded commutative algebra was introduced in \cite{SX} to give the equivalent dual description of a Hom-Lie algebra.

 \begin{defi}{\rm(\cite{SX})}
A {\bf $(\sigma,\tau)$-differential graded commutative algebra}  consists of a graded commutative   algebra $\frkA=\oplus_k\frkA_k$, degree $0$ algebra endomorphisms $\sigma$ and $\tau$, and a degree $1$  operator $d_\frkA:\frkA_k\longrightarrow\frkA_{k+1}$, such that the following conditions are satisfied:
\begin{itemize}
  \item[\rm(1)] $d_\frkA^2=0;$
   \item[\rm(2)] $d_\frkA\circ \sigma=\sigma\circ d_\frkA,\quad d_\frkA\circ \tau=\tau\circ d_\frkA;$
    \item[\rm(3)] $d_\frkA(ab)=d_\frkA(a)\tau(b)+(-1)^{k}\sigma(a)d_\frkA(b),\quad \forall~a\in \frkA_{k},b\in\frkA_{l}.$
\end{itemize}
\end{defi}

 Let $\frkA$ be a commutative algebra and $\sigma$ an algebra isomorphism. Define a skew-symmetric bilinear map $[\cdot,\cdot]_\sigma:\wedge^2\Der_{\sigma,\sigma}(\frkA)\longrightarrow \Der_{\sigma,\sigma}(\frkA)$ by
$$
  [D_1,D_2]_\sigma=\sigma\circ D_1\circ \sigma^{-1}\circ D_2\circ \sigma^{-1}-\sigma\circ D_2\circ \sigma^{-1}\circ D_1\circ \sigma^{-1},\quad \forall D_1,D_2\in\Der_{\sigma,\sigma}(\frkA).
$$
Define $\Ad_\sigma:\Der_{\sigma,\sigma}(\frkA)\rightarrow\Der_{\sigma,\sigma}(\frkA)$  by
$$
  \Ad_\sigma(D)=\sigma\circ D\circ\sigma^{-1}.
$$

\begin{thm}{\rm(\cite{CS-bialgebra})}\label{thm:derhomliealg}
With the above notations,  $(\Der_{\sigma,\sigma}(\frkA),[\cdot,\cdot]_\sigma,\Ad_\sigma)$ is a Hom-Lie algebra.
\end{thm}

  A Hom-Lie algebra $(\g,[\cdot,\cdot]_\g,\phi_\g)$ is called a {\bf regular Hom-Lie algebra} if $\phi_\g$ is
an algebra automorphism. In the sequel, all Hom-Lie algebras are regular.

\begin{defi}{\rm(\cite{SH-rep})}\label{repn}
  A representation of a Hom-Lie algebra $(\g,[\cdot,\cdot]_\g,\phi_\g)$ on
  a vector space $V$ with respect to $\beta\in\gl(V)$ is a linear map
  $\rho:\g\longrightarrow \gl(V)$, such that for all
  $x,y\in \g$, the following equalities are satisfied:
  \begin{eqnarray*}
\rho(\phi_\g(x))\circ \beta&=&\beta\circ \rho(x);\\
    \rho([x,y]_\g)\circ
    \beta&=&\rho(\phi_\g(x))\circ\rho(y)-\rho(\phi_\g(y))\circ\rho(x).
  \end{eqnarray*}
\end{defi}
We denote a representation by  $(V,\beta,\rho)$.
See \cite{AEM,MS1,SH-rep} for more details about representations and cohomologies of Hom-Lie algebras and their applications.

 Let $(\g,[\cdot,\cdot]_\g,\phi_\g)$ be a Hom-Lie algebra. The linear map $\phi_\g:\g\longrightarrow\g$ can be extended to a linear map from $\wedge^k\g\longrightarrow\wedge^k\g$, for which we use the same notation $\phi_\g$ via
$$
\phi_\g(x_1\wedge\cdots\wedge x_k)=\phi_\g(x_1)\wedge\cdots\wedge\phi_\g(x_k).
$$  Furthermore, the bracket operation $[\cdot,\cdot]_\g$ can also be extended to $\wedge^\bullet\frkg$
via
$$
\llceil x_1\wedge\cdots \wedge x_m,y_1\wedge \cdots
\wedge y_n\rrceil_\g
=\sum_{i,j}(-1)^{i+j}[x_i,y_j]_\g \wedge
\phi_\g(x_1\wedge\cdots \widehat{x_i}\cdots \wedge
x_m\wedge y_1\wedge\cdots
\widehat{y_j}\cdots\wedge y_n),
$$
for all $x_1\wedge\cdots \wedge x_m\in\wedge^m \g,~y_1\wedge \cdots
\wedge y_n\in\wedge^n \g$. Consequently, $(\oplus_k \wedge^k\g,\wedge,\llceil\cdot,\cdot\rrceil_\g,\phi_\g)$ is a Hom-Gerstenhaber algebra introduced in \cite{LGT}. There is a natural representation of $\g$ on $\wedge^k\g$, which we denote by $\ad$, given by
$$\ad_xY=\llceil x,Y\rrceil_\g,\quad\forall x\in \g,~Y\in\wedge^k\g.$$

Let $(V,\beta,\rho)$ be a representation of a Hom-Lie algebra $(\g,[\cdot,\cdot]_\g,\phi_\g)$. In the sequel, we always assume that $\beta$ is invertible. Define $\rho^*:\g\longrightarrow\gl(V^*)$ as usual by
$$\langle \rho^*(x)(\xi),u\rangle=-\langle\xi,\rho(x)(u)\rangle,\quad\forall x\in\g,u\in V,\xi\in V^*.$$
However, in general $\rho^*$ is not a representation of $\g$ anymore (see \cite{BM} for details). Define $\rho^\star:\g\longrightarrow\gl(V^*)$ by
\begin{equation}\label{eq:new1}
 \rho^\star(x)(\xi):=\rho^*(\phi_\g(x))\big{(}(\beta^{-2})^*(\xi)\big{)},\quad\forall x\in\g,\xi\in V^*.
\end{equation}
More precisely, we have
\begin{eqnarray}\label{eq:new1gen}
\langle\rho^\star(x)(\xi),u\rangle=-\langle\xi,\rho(\phi_\g^{-1}(x))(\beta^{-2}(u))\rangle,\quad\forall x\in\g, u\in V, \xi\in V^*.
\end{eqnarray}
Then $(V^*,(\beta^{-1})^*,\rho^\star)$ is a representation of $(\g,[\cdot,\cdot]_\g,\phi_\g)$, called the dual representation. See \cite{CS-bialgebra} for more details. In particular, $\ad^\star$ which is given by
 \begin{equation}\label{eq:coadjoint}
  \ad^\star_x\xi=\ad^*_{\phi_\g(x)}(\phi_\g^{-2})^*(\xi),\quad \forall x\in\g, \xi\in\g^*,
 \end{equation}
 is called the {\bf coadjoint representation} of $(\g,[\cdot,\cdot]_\g,\phi_\g)$ on $\g^*$ with respect to $(\phi_\g^{-1})^*$.

\begin{defi}{\rm(\cite{CS-bialgebra})}\label{defi:Hom-Liebi}
 Let $(\g,[\cdot,\cdot]_{\g},\phi_\g)$ and $(\g^*,[\cdot,\cdot]_{\g^*},(\phi_\g^{-1})^*)$ be two Hom-Lie algebras.  $(\g,\g^*)$ is called a {\bf purely Hom-Lie bialgebra} if the following compatibility condition holds:
 \begin{eqnarray*}\label{eq:homLiebi}
 \Delta([x,y]_{\g})=\ad_{\phi_\g^{-1}(x)}\Delta (y)-\ad_{\phi_\g^{-1}(y)}\Delta (x),
 \end{eqnarray*}
 where $\Delta:\g\longrightarrow \wedge^2\g$ is the dual of the Hom-Lie algebra structure $[\cdot,\cdot]_{\g^*}:\wedge^2\g^*\longrightarrow \g^*$ on $\g^*$.
 \end{defi}

 \begin{thm}{\rm(\cite{CS-bialgebra})}\label{thm:doublealg}
 Let $(\g,\g^*)$ be a purely Hom-Lie bialgebra. Then $\Big{(}\g\oplus \g^*,[\cdot,\cdot]_{\bowtie},\phi_\g\oplus(\phi_\g^{-1})^*,(\cdot,\cdot)_+\Big{)}$ is a quadratic Hom-Lie algebra, where the Hom-Lie bracket $[\cdot,\cdot]_{\bowtie}$ and the symmetric bilinear form $(\cdot,\cdot)_+$ are given by
 \begin{eqnarray*}
\label{eq:bracketdouble}[x+\xi,y+\eta]_{\bowtie}
&=&\big{(}[x,y]_\g+\AD_{\xi}^\star y-\AD_{\eta}^\star x\big{)}
+\big{(}[\xi,\eta]_{\g^*}+\ad_{x}^\star\eta-\ad_{y}^\star\xi\big{)},\\
\label{eq:bilinear form}
  (x+\xi,y+\eta)_+&=&\xi(y)+\eta(x), \quad \forall x,y\in\g, \xi,\eta\in\g^*,
\end{eqnarray*}
where $\AD^\star$ is the coadjoint representation of the Hom-Lie algebra  $(\g^*,[\cdot,\cdot]_{\g^*},(\phi_\g^{-1})^*)$.
\end{thm}

\section{Hom-Lie algebroids}

 Let $M$ be a differential manifold and $\varphi:M\longrightarrow M$ a smooth map. Then the pullback map $\aap:C^\infty(M)\longrightarrow C^\infty(M)$ is a morphism of the function ring $C^\infty(M)$, i.e.
 \begin{equation}
   \label{f1}\aap(fg)=\aap(f)\aap(g),\quad \forall f,g\in C^\infty(M).
 \end{equation}

 \begin{defi}
   A {\bf Hom-bundle} is a vector bundle $\aab\rightarrow M$ equipped with a smooth map $\varphi:M\rightarrow M$ and a linear map $\aaa:\aar\rightarrow\aar$ such that the following condition holds:
   \begin{equation}\label{f555}
    \aaa(fx)=\aap(f)\aaa(x),\quad \forall x\in\aar,f\in C^{\infty}(M).
   \end{equation}
  We will call the condition \eqref{f555} {\bf $\varphi^*$-function linear} in the sequel.  A Hom-bundle $(\aab\rightarrow M,\varphi,\aaa)$ is said to be {\bf invertible} if $\varphi$ is a diffeomorphism and $\aaa$ is an invertible linear map.
 \end{defi}

 \begin{ex}{\rm
Let $(\aab\rightarrow M,\varphi,\aaa)$ be a Hom-bundle.  $\aaa$ induces a linear map from $\aaf$ to $\aaf$, which we use the same notation, by
\begin{equation}\label{eq:auto}
\aaa(X)=\aaa(x_1)\wedge\cdots\wedge\aaa(x_k),~\forall~ X=x_1\wedge\cdots\wedge x_k\in\aaf.
 \end{equation}
 Then $(\wedge^kA,\varphi,\aaa)$ is a Hom-bundle.

Assume that $(\aab\rightarrow M,\varphi,\aaa)$ is invertible. We use $\varphi^{-1}$ and $\aaa^{-1}$ to denote the inverses of $\varphi$ and $\aaa$ respectively.
 It is straightforward to obtain
 \begin{equation}
   \label{f4}\aaa^{-1}(fX)=(\aap)^{-1}(f)\aaa^{-1}(X),\quad\forall~f\in C^\infty(M),~ X\in\aaf.
 \end{equation}
 Therefore, $(\aab\rightarrow M,\varphi^{-1},\aaa^{-1})$ is also a Hom-bundle.

 Define $\aad:\aag\rightarrow\aag$ by\footnote{Since $\aaa$ is not a bundle map, the usual $\aaa^*$ is meaningless. }
\begin{eqnarray}\label{eq:auto2}
 (\aad(\Xi))(X)=\aap\Xi(\aaa^{-1}(X)),\quad \forall X\in\aaf, \Xi\in\aag.
 \end{eqnarray}
 By \eqref{f4}, $\aad$ is well-defined. Obviously, we have
 \begin{equation}
   \label{f5}\aad(f\Xi)=\aap(f)\aad(\Xi).
 \end{equation}
 Therefore, $(\wedge^kA^*,\varphi,\aad)$ is a Hom-bundle.
 }
 \end{ex}

  Let $B\longrightarrow M$ be a vector bundle over $M$ and $\varphi:M\longrightarrow M$ a smooth map. We use $\varphi^!B$ to denote the pullback bundle of $B$ along $\varphi$. For any $y\in \Gamma(B)$, we use $y^!\in\Gamma(\varphi^!B)$ to denote the corresponding pullback section, i.e. $y^!(m)=y(\varphi(m))$ for $m\in M$.

  \begin{ex}{\rm(\cite{LGT})
    For any bundle map $\Phi:\varphi^!B\longrightarrow B$, define $\phi_{\varphi^!B} :\Gamma(\varphi^!B)\longrightarrow \Gamma(\varphi^!B)$ by
    $$
    \phi_{\varphi^!B}(x)=(\Phi(x))^!,\quad\forall x\in\Gamma(\varphi^!B).
    $$
    Then we have
    $$
    \phi_{\varphi^!B}(fx)=(\Phi(fx))^!=(f\Phi(x))^!=\varphi^*(f)(\Phi(x))^!=\varphi^*(f)\phi_{\varphi^!B}(x).
    $$
    Therefore, $(\varphi^!B,\varphi,\phi_{\varphi^!B})$ is a Hom-bundle.
    }
  \end{ex}

\begin{ex}{\rm
Let $M$ be a manifold and $\varphi:M\rightarrow M$ a diffeomorphism. As pointed in \cite[Remark 3.8]{LGT},  $\Gamma(\varphi^!TM)$ can be identified with the set of $(\aap,\aap)$-derivations $\Der_{\aap,\aap}(C^{\infty}(M))$ on $C^\infty(M)$, i.e. for all $f,g\in C^{\infty}(M),x\in\Gamma(\varphi^!TM),$ we have
 \begin{eqnarray*}\label{f00}
 x(fg)=x(f)\aap(g)+\aap(f)x(g).
 \end{eqnarray*}
 Define $\Ad_{\aap}:\Gamma(\varphi^!TM)\longrightarrow \Gamma(\varphi^!TM)$ by
 \begin{equation}\label{eq:Auto1}
 \Ad_{\aap}(x)=\aap\circ x\circ(\aap)^{-1},~\forall x\in\Gamma(\varphi^!TM).
 \end{equation}
 Then $(\varphi^!TM,\varphi,\Ad_{\aap} )$ is a Hom-bundle.
 }
\end{ex}

\subsection{Dual description of  Hom-Lie algebroids }

 \begin{defi}\label{defi:hla}
 A {\bf Hom-Lie algebroid} structure on a Hom-bundle $(\aab\rightarrow M,\varphi,\aaa)$ is a pair that consists of a Hom-Lie algebra structure $(\aar,[\cdot,\cdot]_{\aab},\phi_A)$ on the section space $\Gamma(A)$ and a bundle map $a_A:\aab\rightarrow \varphi^!TM$, called the anchor, such that the following conditions are satisfied:
 \begin{itemize}
 \item[\rm(i)]For all $x, y\in\aar$ and $f\in C^{\infty}(M)$, $[x,fy]_\aab=\aap(f)[x,y]_\aab+a_A(\aaa (x))(f)\aaa(y)$;
 \item[\rm(ii)]the anchor $a_A:\aar\longrightarrow \Gamma(\varphi^!TM)$ is a representation of Hom-Lie algebra $(\aar,[\cdot,\cdot]_\aab,\aaa)$ on $C^{\infty}(M)$ with respect to $\aap$.
 \end{itemize}
 We denote a Hom-Lie algebroid by $(\aab,\varphi,\phi_A,[\cdot,\cdot]_A,a_A)$.
 \end{defi}

 In particular, if $\varphi=\idd$ and $\phi_A=\idd$, a Hom-Lie algebroid is exactly a Lie algebroid. See \cite{General theory of Lie groupoid and Lie algebroid} for more details about Lie algebroids.

 \begin{rmk}In  \cite{LGT}, the authors had already defined a Hom-Lie algebroid.
 There is a slight difference  between the above definition of a Hom-Lie algebroid and that one. In \cite{LGT}, a Hom-Lie algebroid is defined to be a quintuple $(\mathscr{A}\rightarrow M, \varphi,\alpha,[\cdot,\cdot]_{\mathscr{A}},\rho)$, where $(\mathscr{A}\rightarrow M, \varphi,\alpha)$ is a Hom-bundle, $(\Gamma(\mathscr{A}),[\cdot,\cdot]_{\mathscr{A}},\alpha)$ is a Hom-Lie algebra, $\rho: \varphi^{!}\mathscr{A}\rightarrow\varphi^!TM $ is a bundle map such that $(C^\infty(M),\varphi^*,\rho)$ is a representation of $(\Gamma(\mathscr{A}),[\cdot,\cdot]_{\mathscr{A}},\alpha)$ and the following condition holds:
 \begin{equation}\label{eq:conold}
   [X,FY]_{\mathscr{A}}=\aap(f)[X,Y]_{\mathscr{A}}+\rho(X)[f]\alpha(Y),\quad \forall X,Y\in\Gamma(\mathscr{A}), F\in C^\infty(M),
 \end{equation}
 where the value of $\rho(X)[F]$ at $m\in M$ is equal to $\langle d_{\varphi(m)}F,\rho_m(X_{\varphi(m)})\rangle$. Note that the condition \eqref{eq:conold} is not the same as the condition ${\rm(i)}$ in Definition \ref{defi:hla}.

  If $\varphi$ is a diffeomorphism,  then $(\aab=\varphi^{!}\mathscr{A},\varphi,\phi_A=\alpha^!,[\cdot,\cdot]_A,a_A=\rho)$ is a Hom-Lie algebroid in the sense of Definition \ref{defi:hla}, where $\alpha^!:\Gamma(\mathscr{A})\longrightarrow \Gamma(\mathscr{A})$ is defined by $\alpha^!(Y^!)=\alpha(Y)^!$ for $Y\in\Gamma(\mathscr{A})$ and $[\cdot,\cdot]_A$ is defined by $[X^!,Y^!]_A=[X,Y]_{\mathscr{A}}^!$ for $X,Y\in\Gamma(\mathscr{A})$. To see that the condition ${\rm(i)}$ holds, we let $x=X^!, y=Y^!\in \Gamma(A)$ and $f=\varphi^*(F)\in C^\infty(M)$. Then we have
  \begin{eqnarray*}
   ~ [x,fy]_A
   &=&[X^!,(FY)^!]_A
   =[X,FY]_{\mathscr{A}}^!\\
    &=&\big(\varphi^*(F)[X,Y]_{\mathscr{A}}\big)^!+\big(\rho(X)[F]\alpha(Y)\big)^!\\
    &=&\varphi^*(\varphi^*(F))[X,Y]_{\mathscr{A}}^!+\varphi^*(\rho(X)[F])(\alpha(Y))^!\\
    &=&\varphi^*(f)[X^!,Y^!]_A+\rho(\alpha(X))[\varphi^*(F)](\alpha^!(Y^!))\\
    &=&\varphi^*(f)[x,y]_A+a_A(\phi_A(x))(f)(\phi_A(y)).
  \end{eqnarray*}
  In the sequel, we will see that our  Hom-Lie algebroids are easy to be treated and have interesting examples.
 \end{rmk}

 In the sequel, we always assume that the underlying Hom-bundle $(\aab\rightarrow M,\varphi,\aaa)$ is invertible.
Define a differential operator
$
 \dM:\aag\rightarrow\Gamma(\wedge^{k+1}A^*)\nonumber
$
by
 \begin{eqnarray}\label{D1}
 &&(\dM\Xi)(x_1,\cdots,x_{k+1})=\sum_{i=1}^{k+1}(-1)^{i+1}a_A(x_i)\Xi\big(\aaa^{-1}(x_1),\cdots,\widehat{\aaa^{-1}(x_i)},\cdots,\aaa^{-1}(x_{k+1})\big)\nonumber\\
 &&\quad+\sum_{i< j}(-1)^{i+j}\aad(\Xi)\big([\aaa^{-1}(x_i),\aaa^{-1}(x_j)]_{\aab},x_1,\cdots,\widehat{x_i},\cdots,\widehat{x_j},\cdots,x_{k+1}\big).
 \end{eqnarray}

 \begin{thm}\label{thm:DGCA1}
 Let   $(\aab,\varphi,\aaa,[\cdot,\cdot]_{\aab},a_A)$ be a Hom-Lie algebroid. Then we have $\dM^2=0.$ Furthermore, for all $\Xi\in\aag$ and $\Theta\in\Gamma(\wedge^l\aab^*)$, there holds:
\begin{eqnarray}
\label{d0}\dM\circ \aad&=&\aad\circ \dM;\\
\label{d2}\dM(\Xi\wedge \Theta)&=&\dM \Xi\wedge\aad(\Theta)+(-1)^k\aad(\Xi)\wedge\dM\Theta.
\end{eqnarray}
Therefore, $(\oplus_k\aag,\wedge,\aad,\dM)$ is a $(\aad,\aad)$-differential graded commutative algebra.
\end{thm}
\pf  Since $a_A$ is a representation of   $(\aar,[\cdot,\cdot]_\aab,\aaa)$ on $C^{\infty}(M)$ with respect to $\aap$, we can obtain $\dM^2=0$ directly.   By direct calculation, \eqref{d0} follows immediately.

In the following, we prove that \eqref{d2} holds by induction on $k$. For $k=0$, for all $f\in C^\infty(M)$, we have
 \begin{eqnarray*}
 &&\dM(f\Theta)(x_1,\cdots,x_{l+1})\nonumber\\
 &=&\sum_{i=1}^{l+1}(-1)^{i+1}a_A(x_i)\big{(}f\Theta(\aaa^{-1}(x_1),\cdots,\widehat{\aaa^{-1}(x_i)},\cdots,\aaa^{-1}(x_{l+1}))\big{)}\nonumber\\
 &&+\sum_{i< j}(-1)^{i+j}(\aad(f\Theta))([\aaa^{-1}(x_i),\aaa^{-1}(x_j)]_{\aab},x_1,\cdots,\widehat{x_i},\cdots,\widehat{x_j},\cdots,x_{l+1})\nonumber\\
 &=&\sum_{i=1}^{l+1}(-1)^{i+1}\big{(}a_A(x_i)(f)\aap\Theta(\aaa^{-1}(x_1),\cdots,\widehat{\aaa^{-1}(x_i)},\cdots,\aaa^{-1}(x_{l+1}))\nonumber\\
 &&+\aap(f)a_A(x_i)(\Theta(\aaa^{-1}(x_1),\cdots,\widehat{\aaa^{-1}(x_i)},\cdots,\aaa^{-1}(x_{l+1})){)}\big)\nonumber\\
 &&+\sum_{i< j}(-1)^{i+j}\aap(f)\aad(\Theta)([\aaa^{-1}(x_i),\aaa^{-1}(x_j)]_{\aab},x_1,\cdots,\widehat{x_i},\cdots,\widehat{x_j},\cdots,x_{l+1})
  \end{eqnarray*}
   \begin{eqnarray*}
 &=&\sum_{i=1}^{l+1}(-1)^{i+1}\dM f(x_i)\aad(\Theta)(x_1,\cdots,\widehat{x_i},\cdots,x_{k+1})+\aap(f)\dM\Theta(x_1,\cdots,x_{l+1})\nonumber\\
 &=&\big{(}\dM f\wedge\aad(\Theta)+\aap(f)\dM\Theta\big{)}(x_1,\cdots,x_{l+1}),\nonumber
 \end{eqnarray*}
 which means that
 \begin{eqnarray*}
 \dM (f\Theta)=\dM f\wedge\aad(\Theta)+\aap(f)\dM\Theta.
 \end{eqnarray*}
Assume that for $k\leq n$, \eqref{d2} holds. For $\xi\in\aaq,\Xi\in\Gamma(\wedge^n A^*)$, we have
\begin{eqnarray*}
\dM(\Xi\wedge\xi\wedge\Theta)
&=&\dM\Xi\wedge\aad(\xi\wedge\Theta)+(-1)^n\aad(\Xi)\wedge\dM(\xi\wedge\Theta)\\
&=&\dM\Xi\wedge(\aad(\xi)\wedge\aad(\Theta))+(-1)^n\aad(\Xi)\wedge(\dM\xi\wedge\aad(\Theta)-\aad(\xi)\wedge\dM\Theta)\\
&=&\Big{(}\dM\Xi\wedge\aad(\xi)+(-1)^n\aad(\Xi)\wedge\dM\xi\Big{)}\wedge\aad(\Theta)
+(-1)^{n+1}\aad(\Xi)\wedge\aad(\xi)\wedge\dM\Theta\\
&=&\dM(\Xi\wedge\xi)\wedge\aad(\Theta)+(-1)^{n+1}\aad(\Xi\wedge\xi)\wedge\dM\Theta,
\end{eqnarray*}
which implies that, for $k=n+1$, \eqref{d2} holds.\qed\vspace{3mm}

 The converse of the above result is also true.

\begin{thm}\label{thm:dgca}
Let $(\aab\rightarrow M,\varphi,\aaa)$ be an invertible Hom-bundle. If $(\oplus_k\aag,\wedge,\aad,\dM)$ is a $(\aad,\aad)$-differential graded commutative algebra, then $(\aab,\varphi,\aaa,[\cdot,\cdot]_{\aab},a_A)$ is a Hom-Lie algebroid, where for all $\xi\in\aaq,x_1,x_2\in\aar$ and $f\in C^\infty(M)$, the anchor $a_A$ and the bracket $[\cdot,\cdot]_\aab$  are respectively defined by
\begin{eqnarray}
\label{eq:DGCA1}a_A(x)(f)&=&\langle\dM f,x\rangle
\end{eqnarray}
and
\begin{eqnarray}
\nonumber\langle[x_1,x_2]_\aab,\xi\rangle&=&-(\dM\aadd(\xi))(\aaa (x_1),\aaa (x_2))\\
&&+a_A(\aaa(x_1)\langle\aadd(\xi),x_2\rangle
-a_A(\aaa(x_2))\langle\aadd(\xi),x_1\rangle.
\label{eq:DGCA2}
\end{eqnarray}
\end{thm}
\pf
For all $x\in\aar, f\in C^{\infty}(M)$, by \eqref{eq:DGCA1}, we have
\begin{eqnarray*}
a_A(\aaa (x))(\varphi^*(f))
&=&\langle\dM \varphi^*(f),\aaa (x)\rangle=\langle\aad(\dM f),\aaa (x)\rangle=\varphi^*\langle\dM f,x\rangle=\varphi^*a_A(x)(f),
\end{eqnarray*}
which implies that
\begin{equation}\label{eq:DGCA3}
a_A(\aaa (x))\circ\aap=\aap\circ a_A(x).
\end{equation}

For all  $x_1,x_2\in\Gamma(A)$ and $f\in\CWM$, we have
$$(\dM^2 f)(\aaa (x_1),\aaa (x_2))=0,$$
which implies that
\begin{equation}\label{eq:DGCA4}
a_A([x_1,x_2]_\aab)\circ\varphi^*=a_A(\aaa (x_1))\circ a_A(x_2)-a_A(\aaa (x_2))\circ a_A(x_1).
\end{equation}

For all $\xi\in\aaq,x_1,x_2\in\aar$, by \eqref{eq:DGCA2}, we have
\begin{eqnarray*}
\dM(\aad(\xi))(x_1,x_2)&=&a_A(x_1)(\aad(\xi)(\aaa^{-1}(x_2)))-a_A(x_2)(\aad(\xi)(\aaa^{-1}(x_1)))\\
&&-((\aad)^2(\xi))([\aaa^{-1}(x_1),\aaa^{-1}(x_2)]_\aab);\\
\aad\dM\xi(x_1,x_2)&=&\aap\dM\xi(\aaa^{-1}(x_1),\aaa^{-1}(x_2))\\
&=&a_A(x_1)(\aad(\xi)(\aaa^{-1}(x_2)))-a_A(x_2)(\aad(\aaa^{-1}(x_1)))\\
&&-((\aad)^2(\xi))(\aaa[\aaa^{-2}(x_1),\aaa^{-2}(x_2)]_A).
\end{eqnarray*}
Then by \eqref{d0}, we have
\begin{equation}\label{eq:DGCA5}
\aaa[x_1,x_2]_\aab=[\aaa(x_1),\aaa(x_2)]_\aab,
\end{equation}
which implies that $\aaa$ is an algebra endomorphism.

By \eqref{d2} and  \eqref{eq:DGCA2},
for all $\Xi\in\Gamma(\wedge^2\aab^*),$ the following equality holds:
\begin{eqnarray*}
\dM(\Xi)(x_1,x_2,x_3)&=&a_A(x_1)(\Xi(\aaa^{-1}(x_2),\aaa^{-1}(x_3)))-a_A(x_2)(\Xi(\aaa^{-1}(x_1),\aaa^{-1}(x_3)))\\
&&+a_A(x_3)(\Xi(\aaa^{-1}(x_1),\aaa^{-1}(x_2)))-\aad(\Xi)([\aaa^{-1}(x_1),\aaa^{-1}(x_2)]_\aab,x_3)\\
&&+\aad(\Xi)([\aaa^{-1}(x_1),\aaa^{-1}(x_3)]_\aab,x_2)-\aad(\Xi)([\aaa^{-1}(x_2),\aaa^{-1}(x_3)]_\aab,x_1).
\end{eqnarray*}
 Then by  $(\dM^2\xi)(x_1,x_2,x_3)=0$, for all $\xi\in\aaq$ and $x_1,x_2,x_3\in \aar$, we have
\begin{eqnarray}\label{eq:DGCA6}
~[[x_1,x_2]_\aab,\aaa (x_3)]_\aab+[[x_2,x_3]_\aab,\aaa (x_1)]_\aab+[[x_3,x_1]_\aab,\aaa (x_2)]_\aab=0.
\end{eqnarray}
Thus, by \eqref{eq:DGCA5} and \eqref{eq:DGCA6}, $(\aar,[\cdot,\cdot]_\aab,\aaa)$ is a Hom-Lie algebra. By \eqref{eq:DGCA3} and \eqref{eq:DGCA4}, $a_A$ is a representation of Hom-Lie algebra $(\aar,[\cdot,\cdot]_\aab,\aaa)$ on $C^{\infty}(M)$ with respect to $\aap$.

For $\xi\in\Gamma(A^*),x_1,x_2\in\Gamma(A)$ and $f\in\CWM$, we have
\begin{eqnarray*}
\dM\xi(\phi_A(x_1),\phi_A(fx_2))&=&\varphi^*(f)a_A(\phi_A(x_1))(\xi(x_2))+a_A(\phi_A(x_1))(f)\aap(\xi(x_2))\\
&&-\varphi^*(f)a_A(\phi_A(x_2))(\xi(x_1))-\aad(\xi)([x_1,fx_2]_\aab)\\
&=&\varphi^*(f)a_A(\phi_A(x_1))(\xi(x_2))-\varphi^*(f)a_A(\phi_A(x_2))(\xi(x_1))\\
&&+a_A(\phi_A(x_1))(f)\aad(\xi)(\aaa(x_2))-\aad(\xi)([x_1,fx_2]_\aab).
\end{eqnarray*}
On the other hand, we have
\begin{eqnarray*}
&&\dM\xi(\phi_A(x_1),\phi_A(fx_2))=\varphi^*(f)\dM\xi(\phi_A(x_1),\phi_A(x_2))\\
&=&\varphi^*(f)a_A(\phi_A(x_1))(\xi(x_2))-\varphi^*(f)a_A(\phi_A(x_2))(\xi(x_1))-\varphi^*(f)\aad(\xi)([x_1,x_2]_\aab).
\end{eqnarray*}
Thus, we have
\begin{eqnarray*}
~[x_1,fx_2]_\aab=\aap(f)[x_1,x_2]_\aab+a_A(\aaa (x_1))(f)\aaa (x_2).
\end{eqnarray*}

Therefore, $(\aab,\varphi,\aaa,[\cdot,\cdot]_{\aab},a_A)$ is a Hom-Lie algebroid.\qed

\subsection{Tangent Hom-Lie algebroids and action Hom-Lie algebroids }

 Let $\varphi:M\rightarrow M$ be a diffeomorphism. Define a skew-symmetric bilinear operation $[\cdot,\cdot]_{\varphi^*}:\wedge^2\Gamma(\varphi^!TM)\longrightarrow \Gamma(\varphi^!TM)$ by
 \begin{equation}\label{brac}
   [x,y]_{\varphi^*}=\aap\circ x\circ(\aap)^{-1}\circ y\circ (\aap)^{-1}-\aap\circ y\circ(\aap)^{-1}\circ x\circ (\aap)^{-1}.
 \end{equation}

 \begin{pro}\label{mmm}
With the above notations, $(\varphi^!TM,\varphi,\Ad_{\aap},[\cdot,\cdot]_{\varphi^*}, \idd)$ is a Hom-Lie algebroid, where $\Ad_{\aap}$ is given by \eqref{eq:Auto1} and $[\cdot,\cdot]_{\varphi^*}$ is given by \eqref{brac}.
 \end{pro}
\pf By Theorem \ref{thm:derhomliealg}, $(\Gamma(\varphi^!TM),[\cdot,\cdot]_{\varphi^*},\Ad_{\aap})$ is a Hom-Lie algebra. For all $f\in\CWM$, we have
 \begin{eqnarray*}
 ~[x,fy]_{\varphi^*}(g)&=&\big{(}\aap\circ x\circ(\aap)^{-1}\circ fy\circ(\aap)^{-1}
 -\aap\circ fy\circ(\aap)^{-1}\circ x\circ(\aap)^{-1}\big{)}(g)\\
 &=&\aap x\big{(}(\aap)^{-1}(f)(\aap)^{-1}y((\aap)^{-1}g)\big{)}
 -\aap(f)\aap y((\aap)^{-1}x((\aap)^{-1}g))\\
 &=&\aap x((\aap)^{-1}(f))\aap y((\aap)^{-1}g)
 +\aap(f)\aap x((\aap)^{-1}y((\aap)^{-1}g))\\
&&-\aap(f)\aap y((\aap)^{-1}x((\aap)^{-1}g))\\
&=&\big{(}\aap(f)[x,y]_{\varphi^*}+\Ad_{\aap}(x)(f)\Ad_{\aap}(y)\big{)}(g),
 \end{eqnarray*}
 which implies that $[x,fy]_{\varphi^*}=\aap(f)[x,y]_{\varphi^*}+\Ad_{\aap}(x)(f)\Ad_{\aap}(y).$

It is obvious that $\idd$ is a representation of $(\Gamma(\varphi^!TM),[\cdot,\cdot]_{\varphi^*},\Ad_\aap)$ on $C^{\infty}(M)$ with respect to $\aap.$ Thus $(\varphi^!TM,\varphi,\Ad_{\aap},[\cdot,\cdot]_{\varphi^*}, \idd)$ is a Hom-Lie algebroid.\qed

 \begin{rmk}
 The above proposition tells us that the pullback of the tangent Lie algebroid $TM$ is a Hom-Lie algebroid. We can obtain a more general result: for any Lie algebroid equipped with some extra structures, we can obtain a Hom-Lie algebroid.
  Let $(\huaA,[\cdot,\cdot]_{\huaA},a_\huaA)$ be a Lie algebroid and $(\huaA,\varphi,\alpha)$ a Hom-bundle. If $\alpha[x,y]_{\huaA}=[\alpha(x),\alpha(y)]_{\huaA}$ and $a_\huaA(\alpha(x))\circ\varphi^*=\varphi^*\circ a_\huaA(x)$, then $(\aab=\varphi^{!}\huaA,\varphi,\phi_A=\alpha^!,[\cdot,\cdot]_A,a_A)$ is a Hom-Lie algebroid, where $\alpha^!:\Gamma({A})\longrightarrow \Gamma({A})$,  the Hom-Lie bracket $[\cdot,\cdot]_A$ and $a_A:\Gamma({A})\longrightarrow\Gamma(\varphi^!TM)$ are defined by $$\alpha^!(x^!)=\alpha(x)^!,\quad [x^!,y^!]_A=[\alpha(x),\alpha(y)]_{\huaA}^!,\quad a_A(x^!)=\varphi^*\circ a_\huaA(x)=a_\huaA(x)^!$$
  for all $x,y\in\Gamma(\huaA)$. In particular, let $(\huaA,[\cdot,\cdot]_{\huaA},a_\huaA)$ be the tangent Lie algebroid $TM$ and $\alpha=\Ad_{\varphi^*}$, we obtain the Hom-Lie algebroid $(\varphi^!TM,\varphi,\Ad_{\aap},[\cdot,\cdot]_{\varphi^*}, \idd)$ given in the above proposition.
 \end{rmk}

\begin{defi}
Let $(A,\varphi,\phi_A,[\cdot,\cdot]_A,a_A)$ and $(B,\varphi,\phi_B,[\cdot,\cdot]_B,a_B)$ be two Hom-Lie algebroids over the same base $M$. A bundle map $\sigma:A\longrightarrow B$ is called a {\bf morphism} of Hom-Lie algebroids, if for all $x,y\in\Gamma(A)$, the following conditions are satisfied:
\begin{eqnarray*}
a_B\circ \sigma&=&a_A,\\
\sigma\circ\phi_A(x)&=&\phi_B\circ\sigma(x),\\
\sigma([x,y]_A)&=&[\sigma(x),\sigma(y)]_B.
\end{eqnarray*}
\end{defi}

 \begin{pro}
  Let $(\aab,\varphi,\aaa,[\cdot,\cdot]_{\aab},a_A)$ be a Hom-Lie algebroid and $\varphi:M\rightarrow M$ be a diffeomorphism. Then the anchor $a_A$ is a morphism from the Hom-Lie algebroid $(\aab,\varphi,\aaa,[\cdot,\cdot]_{\aab},a_A)$ to the Hom-Lie algebroid $\varphi^!TM$.
 \end{pro}
 \pf Since $a_A:\aar\longrightarrow \Gamma(\varphi^!TM)$ is a representation of $(\aar,[\cdot,\cdot]_\aab,\aaa)$ on $C^{\infty}(M)$ with respect to $\aap$, we have
 $$a_A(\phi_A(x))=\Ad_{\varphi^*}(a_A(x))$$
 and
   \begin{eqnarray*}
   a_A([x,y]_\aab)
   &=&\big{(}a_A(\aaa x)a_A(y)-a_A(\aaa y)a_A(x)\big{)}\circ(\aap)^{-1}\\
   &=&\aap\circ a_A(x)\circ(\aap)^{-1}\circ a_A(y)\circ(\aap)^{-1}
   -\aap\circ a_A(y)\circ(\aap)^{-1}\circ a_A(x)\circ(\aap)^{-1}\\
   &=&[a_A(x),a_A(y)]_{\varphi^*},
   \end{eqnarray*}
 which implies that $a_A$ is a morphism between Hom-Lie algebroids.\qed\vspace{3mm}

 Using the Hom-Lie algebroid $\varphi^!TM$, we can define an action of a Hom-Lie algebra on a manifold.

\begin{defi}
Let $(\frkg,[\cdot,\cdot]_{\frkg},\phi_\g)$ be Hom-Lie algebra and $\varphi:M\rightarrow M$ a diffeomorphism. An {\bf action}
of $\frkg$ on $M$ is a linear map $\rho:\frkg\longrightarrow\Gamma(\varphi^!TM)$, such that for all $x,y\in\frkg$, we have
 \begin{eqnarray*}
 \rho(\phi_\g(x))&=&\Ad_{\aap}(\rho(x));\\
 \rho([x,y]_{\frkg})&=&[\rho(x),\rho(y)]_{\varphi^*}.
 \end{eqnarray*}
\end{defi}

 Given an action of Hom-Lie algebra
$(\frkg,[\cdot,\cdot]_\g,\phi_\g)$ on $M$, let $A=M\times
\g$ be the trivial bundle. Define a linear map $\phi_A:\Gamma(A)\longrightarrow \Gamma(A)$, an anchor map $a_\rho:A\longrightarrow \varphi^!TM$ and a bracket $[\cdot,\cdot]_\rho:\Gamma(A)\times \Gamma(A)\longrightarrow \Gamma(A)$  by
\begin{eqnarray*}
\phi_A(fc_x)&=&\varphi^*(f)c_{\phi_\g(x)},\quad \forall~ x\in\g,f\in\CWM,\label{action hom1}\\
a_\rho(m,x)&=&\rho(x)_m,\quad \forall ~m\in M, x\in\g,\label{action hom2}\\
{[fc_{x},gc_{y}]}_\rho&=&\varphi^*(fg)c_{[x,y]_{\frkg}}+\varphi^*(f)a_\rho(\phi_\g(x))(g)c_{\phi_\g(y)}-\varphi^*(f)a_\rho(\phi_\g(y))(f)c_{\phi_\g(x)},\label{action hom3}
\end{eqnarray*}
where $c_x,c_y$ denote the constant sections of $M\times\frkg$ given by $m\mapsto (m,x)$ and $m\mapsto (m,y)$ respectively.

By straightforward computations, we obtain
\begin{pro}
With the above notations, $(A=M\times\frkg,\varphi,\phi_A,[\cdot,\cdot]_\rho,a_\rho)$ is a
Hom-Lie algebroid, we call it an {\bf action Hom-Lie algebroid}.
\end{pro}
\subsection{Differential calculus on Hom-Lie algebroids }
In this subsection, we give the Lie derivative for a Hom-Lie algebroid for later applications.

Given a Hom-Lie algebroid $(\aab,\varphi,\aaa,[\cdot,\cdot]_{\aab},a_A)$, define the bilinear operation
$\llceil\cdot,\cdot\rrceil_A:\Gamma(\wedge^pA)\times\Gamma(\wedge^qA)\longrightarrow\Gamma(\wedge^{p+q-1}A)$
by
\begin{eqnarray}
\llceil x_1\wedge\cdots\wedge x_p,y_1\wedge\cdots\wedge y_q\rrceil_A
&=&\sum_{i=1}^p\sum_{j=1}^q(-1)^{i+j}[x_i,y_j]_A\wedge\phi_A(x_1)\wedge\cdots\wedge\widehat{\phi_A(x_i)}\wedge\cdots\wedge\phi_A(x_p)\nonumber\\
&&\wedge\phi_A(y_1)\wedge\cdots\wedge\widehat{\phi_A(y_j)}\wedge\cdots\wedge\phi_A(y_q)
\end{eqnarray}
and
\begin{eqnarray}
\llceil x,f\rrceil_A=a_A(\phi_A(x))(f).
\end{eqnarray}
Then $(\oplus_k\Gamma(\wedge^kA),\wedge,\llceil\cdot,\cdot\rrceil_A,\phi_A)$ is a Hom-Gerstenhaber-algebra \cite{LGT}.

The bracket $\llceil\cdot,\cdot\rrceil_A$ is called a \textbf{Hom-Schouten bracket}. Actually the Hom-Schouten bracket is determined by the following properties,
\begin{itemize}
\item[$\rm(a)$]$\llceil\Gamma(\wedge^kA),\Gamma(\wedge^lA)\rrceil_A\subset\Gamma(\wedge^{k+l-1}A)$ for all
$k$ and $l$;
\item[$\rm(b)$]for all smooth function $f\in C^\infty(M)$
and $x\in\Gamma(A)$,
\begin{eqnarray*}
{\llceil x, f\rrceil_A}=a_A(\phi_A(x))(f);
\end{eqnarray*}
\item[$\rm(c)$]for all $x,y\in\Gamma(A)$, $\llceil x, y\rrceil_A=[x, y]_A;$
\item[$\rm(d)$]for all $X\in\Gamma(\wedge^kA)$ and $Y\in\Gamma(\wedge^lA)$,
\begin{eqnarray*}
\llceil X,Y\rrceil_A&=&-(-1)^{(k-1)(l-1)}\llceil Y,X\rrceil_A,\\
{\llceil X,Y\wedge Z\rrceil_A}&=&\llceil X,Y\rrceil_A\wedge \phi_A(Z)+(-1)^{(k-1)l}\phi_A(Y)\wedge\llceil Y,Z\rrceil_A.
\end{eqnarray*}
\end{itemize}

For all $X\in\Gamma(\wedge^kA)
$, define the {\bf interior multiplication}
$i_X:\Gamma(\wedge^mA^*)\rightarrow\Gamma(\wedge^{m-k}A^*)$
  by
 \begin{eqnarray}\label{D2}
 (i_X\Xi)(y_1,\cdots,y_{m-k})=(\aad(\Xi))(\aaa(X),y_1,\cdots,y_{m-k}).
 \end{eqnarray}
For all $f\in C^{\infty}(M), X\in\Gamma(\wedge^kA),\Xi\in\Gamma(\wedge^mA^*),$ we have
$$i_{f X}\Xi=i_X(f\Xi)=\aap(f) i_X\Xi.$$

 By \eqref{f5} and \eqref{D2}, we get the following formulas directly.
 \begin{pro}
 For all $X\in\Gamma(\wedge^kA),Y\in\Gamma(\wedge^lA),\Xi\in\Gamma(\wedge^mA^*)$ we have
\begin{eqnarray*}
\label{L11}i_{\aaa(X)}\circ i_Y&=&i_{\aaa(Y\wedge X)}\circ\aad=(-1)^{kl}i_{\aaa (Y)}\circ i_X,\\
\label{L20}\aad(i_X\Xi)&=&i_{\aaa (X)}\aad(\Xi),\\
( {\aad})^{-1}(i_X\Xi)&=&i_{\phi_A^{-1} (X)}{(\aad)}^{-1}(\Xi).
\end{eqnarray*}
\end{pro}

For $x\in\Gamma(A)$ and $X\in\Gamma(\wedge^k A)$, we write
$$L_xX=\llceil x,X\rrceil_A,$$
which is the {\bf Lie derivative of multi-sections}. The Lie derivative $L_x:\Gamma(\wedge^k A)\longrightarrow \Gamma(\wedge^k A)$ has the following properties.

\begin{pro}
For all $X\in\Gamma(\wedge^k A),Y\in\Gamma(\wedge^l A)$ and $f\in\CWM$, we have
\begin{eqnarray}
\label{eq:LieD1}L_x(X\wedge Y)&=&L_xX\wedge \phi_A(Y)+\phi_A(X)\wedge L_xY,\\
\label{eq:LieD2}L_{[x,y]_A}\circ \phi_A&=&L_{\phi_A(x)}L_y-L_{\phi_A(y)}L_x,\\
\label{eq:LieD3}L_x (fX)&=&\varphi^*(f)L_x X+a_A(\phi_A(x))\phi_A(X),\\
\label{eq:LieD4}L_{fx}X&=&\varphi^*(f)L_x X-\phi_A(x)\wedge i_{(\aad)^{-1}(\dM f)}X.
\end{eqnarray}
\end{pro}
\pf We only give the proof of \eqref{eq:LieD4}. Others follow directly from the properties of the Hom-Schouten bracket. For all $x,y_1,y_2,\cdots,y_{n}\in \Gamma(A),~~f\in \CWM$, without loss of generality we can assume that $X=y_1\wedge\cdots\wedge y_n$, then we have
\begin{eqnarray*}
L_{fx}(y_1\wedge\cdots\wedge y_n)&=&\sum_{i=1}^{n}\phi_A(y_1)\wedge\cdots\wedge[fx,y_i]_A\wedge\cdots\wedge\phi_A(y_n)\\
&=&\varphi^*(f)\big(\sum_{i=1}^{n}\phi_A(y_1)\wedge\cdots\wedge[x,y_i]_A\wedge\cdots\wedge\phi_A(y_n)\big)\\
&&-\sum_{i=1}^{n}a_A(\phi_A(y_i))(f)\big(\phi_A(y_1)\wedge\cdots\wedge\phi_A(x)\wedge\cdots\wedge\phi_A(y_n)\big)\\
&=&\varphi^*(f)L_{x}(y_1\wedge\cdots\wedge y_n)-\phi_A(x)\wedge i_{(\aad)^{-1}(\dM f)}(y_1\wedge\cdots\wedge y_n),
\end{eqnarray*}
which implies that \eqref{eq:LieD4} holds.\qed
\vspace{3mm}

For all $X\in\Gamma(\wedge^kA)$, define the {\bf Lie derivative}
$L_X:\Gamma(\wedge^m A^*)\rightarrow\Gamma(\wedge^{m-k+1} A^*)$
  by
  \begin{eqnarray}\label{eq:Lied}
  L_X\circ\aad=i_X\circ\dM-(-1)^{k}\dM\circ i_{\aaa^{-1}(X)}.
  \end{eqnarray}
  In particular, for $x\in\Gamma(A), \Xi\in\Gamma(\wedge^m A^*)$,  we have
 \begin{eqnarray}\label{D3}
 \langle L_x\Xi,Y\rangle=a_A(\aaa(x))\langle\Xi,\aaa^{-1}(Y)\rangle-\langle\aad(\Xi), L_x\aaa^{-1}(Y)\rangle, \quad \forall Y\in\Gamma(\wedge^m A).
 \end{eqnarray}
It is straightforward to obtain the following formulas:
\begin{eqnarray}
\label{eq:Lied7}L_{fX}\Xi&=&\aap(f) L_X\Xi-(-1)^k\dM f\wedge i_X\Xi,\quad \forall f\in \CWM, X\in\Gamma(\wedge^kA),~\Xi\in\Gamma(\wedge^m A^*),\\
\label{eq:Lied8}L_x(f\Xi)&=&\aap(f) L_x\Xi+a(\aaa (x))(f)\aad(\Xi),\quad \forall f\in \CWM, x\in\Gamma(A),~\Xi\in\Gamma(\wedge^m A^*).
\end{eqnarray}

\begin{rmk}
When $M$ is a point, the vector bundle $A$ reduces to a vector space. In this case, $\varphi=\idd$ and $a_A=0$. By \eqref{D3}, we have
  \begin{eqnarray*}
  L_x\xi=\ad_{\phi_A(x)}^*(\phi_A^{-2})^*(\xi)=\ad_x^\star\xi,\quad  \forall x\in A, \xi\in A^*.
  \end{eqnarray*}
 Thus, the Lie derivative $L$ for a Hom-Lie algebroid is a natural generalization of the coadjoint representation $\ad^\star$ given by \eqref{eq:coadjoint} for a Hom-Lie algebra.
\end{rmk}

\begin{lem}
For all $x\in\aar,\Xi\in \Gamma(\wedge^mA^*),\Theta\in\Gamma(\wedge^nA^*)$, we have
\begin{eqnarray}
\label{eq:Lied0}i_x(\Xi\wedge\Theta)&=&i_x\Xi\wedge\aad(\Theta)+(-1)^m\aad(\Xi)\wedge i_x\Theta,\\
\label{eq:Lied1}L_x(\Xi\wedge\Theta)&=&L_x\Xi\wedge\aad(\Theta)+\aad(\Xi)\wedge L_x\Theta.
\end{eqnarray}
\end{lem}
\pf By direct calculation, we get \eqref{eq:Lied0}.  By \eqref{d0}, \eqref{d2}
and \eqref{eq:Lied0}, we have
\begin{eqnarray*}
L_x(\Xi\wedge\Theta)&=&i_x\dM(\aad)^{-1}(\Xi\wedge\Theta)+\dM i_{(\aaa)^{-1}(x)}(\aad)^{-1}(\Xi\wedge\Theta)\\
&=&i_x\big(\dM(\aad)^{-1}(\Xi)\wedge\Theta)+(-1)^m\Xi\wedge\dM(\aad)^{-1}(\Theta)\big)\\
&&+\dM\big( i_{(\aaa)^{-1}(x)}(\aad)^{-1}(\Xi)\wedge\Theta+(-1)^m \Xi\wedge i_{(\aaa)^{-1}(x)}(\aad)^{-1}(\Theta) \big)\\
&=&i_x\dM(\aad)^{-1}(\Xi)\wedge\aad(\Theta)+(-1)^{m+1}\aad(\dM(\aad)^{-1}(\Xi))\wedge i_x\Theta\\
&&+(-1)^mi_x\Xi\wedge\aad(\dM(\aad)^{-1}(\Theta))+\aad(\Xi)\wedge i_x\dM(\aad)^{-1}(\Theta)\\
&&+\dM i_{(\aaa)^{-1}(x)}(\aad)^{-1}(\Xi)\wedge\aad(\Theta)-(-1)^m\aad(i_{(\aaa)^{-1}(x)}(\aad)^{-1}(\Xi))\wedge\dM\Theta\\
&&+(-1)^m\dM\Xi\wedge\aad(i_{(\aaa)^{-1}(x)}(\aad)^{-1}(\Theta))+\aad(\Xi)\wedge \dM i_{\phi_A^{-1}(x)}(\aad)^{-1}(\Theta)\\
&=&(i_x\dM(\aad)^{-1}(\Xi)+\dM i_{(\aaa)^{-1}(x)}(\aad)^{-1}(\Xi))\wedge\aad(\Theta)\\
&&+\aad(\Xi)\wedge(i_x\dM(\aad)^{-1}(\Theta)+\dM i_{\phi_A^{-1}(x)}(\aad)^{-1}(\Theta))\\
&=&L_x\Xi\wedge\aad(\Theta)+\aad(\Xi)\wedge L_x\Theta,
\end{eqnarray*}
which implies that \eqref{eq:Lied1} holds.\qed

\begin{lem}
For all $X\in\Gamma(\wedge^k A),Y\in\Gamma(\wedge^l A)$, we have
\begin{equation}
\label{eq:Lied6}i_{\llceil X,Y\rrceil_A}\circ\aad=(-1)^{(k-1)(l-1)}(L_{\aaa (X)}\circ i_Y-(-1)^{(k-1)l}i_{\aaa (Y)}\circ L_X),\\
\end{equation}
\end{lem}
\pf It suffices to show that the bracket defined by \eqref{eq:Lied6} has the same algebraic properties as the Hom-Schouten bracket. We only give the proof of the following equality:
\begin{equation}\label{eq:der1}
{\llceil X,Y\wedge Z\rrceil_A }=\llceil X,Y\rrceil_A \wedge \phi_A(Z)+(-1)^{(k-1)l}\phi_A(Y)\wedge \llceil X,Z\rrceil_A .
\end{equation}
Others can be proved similarly. In fact, for all $Z\in\Gamma(\wedge^mA),\Xi\in\Gamma(A^{k+l+m-1})$, by \eqref{L11}, we have
\begin{eqnarray*}
&&(-1)^{(k-1)(l+m-1)}i_{\llceil X,Y\wedge Z\rrceil_A }\Xi\\
&=&L_{\aaa (X)} i_{Y\wedge Z}(\aad)^{-1}(\Xi)-(-1)^{(k-1)(l+m)}i_{\aaa (Y\wedge Z)} L_X(\aad)^{-1}(\Xi)\\
&=&L_{\aaa (X)} i_{Y\wedge Z}(\aad)^{-1}(\Xi)-(-1)^{(k-1)(l+m)}i_{\aaa (Y\wedge Z)} L_X(\aad)^{-1}(\Xi)\\
&&+(-1)^{(k-1)m}i_{{\aaa (Z)}} L_X i_{(\aaa)^{-1}(Y)}(\aad)^{-2}(\Xi)\\
&&-(-1)^{(k-1)m}i_{{\aaa (Z)}} L_X i_{(\aaa)^{-1}(Y)}(\aad)^{-2}(\Xi)\\
&=&L_{\aaa (X)} i_{Z}(\aad)^{-1}(i_Y(\aad)^{-1}(\Xi))-(-1)^{(k-1)m}i_{{\aaa (Z)}} L_X(\aad)^{-1}(i_Y(\aad)^{-1}(\Xi))\\
&&+(-1)^{(k-1)m}i_{\aaa (Z)}(L_X i_{(\aaa)^{-1}(Y)}(\aad)^{-2}(\Xi)-(-1)^{(k-1)l}i_{Y}L_{(\aaa)^{-1}(X)}(\aad)^{-2}(\Xi))\\
&=&(-1)^{(k-1)(m-1)}i_{\phi_A(Y)\wedge \llceil X,Z\rrceil_A }\Xi+(-1)^{(k-1)(l+m-1)}i_{\llceil X,Y\rrceil_A \wedge \phi_A(Z)}\Xi\\
&=&(-1)^{(k-1)(l+m-1)}i_{(\llceil X,Y\rrceil_A \wedge \phi_A(Z)+(-1)^{(k-1)l}\phi_A(Y)\wedge \llceil X,Z\rrceil_A )}\Xi,
\end{eqnarray*}
which implies that \eqref{eq:der1} holds.\qed

\begin{cor}
For $\pi\in\Gamma(\wedge^2 A),f,g,h\in\CWM$, we have
\begin{eqnarray}\label{eq:HomPFormu}
i_{\llceil\pi,\pi\rrceil_A}(\dM f\wedge\dM g\wedge\dM h)
=-2\Big((\varphi^*)^2\pi\big(\dM(\varphi^*)^{-1}\pi(\dM{(\varphi^*)^{-1}}(f),\dM{(\varphi^*)^{-1}}(g)),\dM{(\varphi^*)^{-1}}(h)\big)+c.p.\Big).
\end{eqnarray}
\end{cor}
\pf
By \eqref{eq:Lied6}, we have
\begin{eqnarray*}
&&i_{\llceil\pi,\pi\rrceil_A}(\dM f\wedge\dM g\wedge\dM h)\\&=&(i_{\aaa (\pi)}\circ L_\pi\circ(\aad)^{-1} -L_{\aaa (\pi)}\circ i_\pi\circ(\aad)^{-1})(\dM f\wedge\dM g\wedge\dM h)\\
&=&i_{\aaa (\pi)}\circ i_\pi\circ\dM(\dM (\varphi^*)^{-2}f\wedge\dM (\varphi^*)^{-2}g\wedge\dM (\varphi^*)^{-2}h))\\
&&-2i_{\aaa (\pi)}\circ\dM\circ i_{\aaa^{-1} (\pi)}(\varphi^*)^{-1}(\dM (\varphi^*)^{-1}f\wedge\dM (\varphi^*)^{-1}g\wedge\dM (\varphi^*)^{-1}h))\\
&&+\dM\circ i_\pi\circ(\varphi^*)^{-1}\circ i_\pi\circ(\varphi^*)^{-1}(\dM f\wedge\dM g\wedge\dM h)\\
&=&-2i_{\aaa (\pi)}\circ\dM\circ i_{\aaa^{-1} (\pi)}(\varphi^*)^{-1}(\dM (\varphi^*)^{-1}f\wedge\dM (\varphi^*)^{-1}g\wedge\dM (\varphi^*)^{-1}h))\\
&=&-2i_{\aaa (\pi)}\circ\dM(\pi(\dM (\varphi^*)^{-1}f,\dM (\varphi^*)^{-1}g)\wedge\dM (\varphi^*)^{-1}h+c.p.)\\
&=&-2(\aad)^2i_{\aaa (\pi)}(\aad)^{-1}\dM(\aad)^{-1}(\pi(\dM (\varphi^*)^{-1}f,\dM (\varphi^*)^{-1}g)\wedge\dM (\varphi^*)^{-1}h+c.p.)\\
&=&-2\Big((\varphi^*)^2\pi\big(\dM(\varphi^*)^{-1}\pi(\dM{(\varphi^*)^{-1}}(f),\dM{(\varphi^*)^{-1}}(g)),\dM{(\varphi^*)^{-1}}(h)\big)+c.p.\Big).\qed
\end{eqnarray*}

\begin{lem}
For all $X\in\Gamma(\wedge^k A),Y\in\Gamma(\wedge^lA),\Xi\in\Gamma(\wedge^mA^*)$, we have
\begin{eqnarray}
\label{eq:Lied2}\aad(L_X\Xi)&=&L_{\aaa (X)}\aad(\Xi),\\
\label{eq:Lied3}L_{\llceil X,Y\rrceil_A}\circ\aad&=&(-1)^{(k-1)(l-1)}(L_{\aaa (X)} L_Y-(-1)^{(k-1)(l-1)}L_{\aaa (Y)} L_X).
\end{eqnarray}
In particular,  $( \aag,\aad, L) $ is a representation of the Hom-Lie algebra $(\Gamma(A),[\cdot,\cdot]_A,\phi_A)$.
\end{lem}
\pf \eqref{eq:Lied2} follows from \eqref{d0} and \eqref{L20}. \eqref{eq:Lied3} follows from \eqref{eq:Lied6}. We omit  details.\qed

\begin{cor}
For all $X\in\Gamma(\wedge^k A)$, we have
\begin{eqnarray}
\label{eq:Lied4}(\aad)^{-1}\circ L_X &=&L_{\phi_A^{-1} (X)}\circ (\aad)^{-1},\\
\label{eq:Lied5}\dM\circ L_X&=&-(-1)^{k}L_{\aaa (X)}\circ\dM.
\end{eqnarray}
\end{cor}

\subsection{The Hom-Lie algebroid associated to a Hom-Poisson manifold}
Recall that $(\varphi^!TM,\varphi,\Ad_{\aap},[\cdot,\cdot]_{\varphi^*},\idd)$ is a Hom-Lie algebroid,  where $\Ad_{\aap}$ and $[\cdot,\cdot]_{\varphi^*}$ are given by \eqref{eq:Auto1} and \eqref{brac} respectively.

\begin{defi}
A bisection $\pi\in\Gamma(\wedge^2 \varphi^!TM)$ on a manifold $M$ is said to be {\bf Hom-Poisson tensor} if $\llceil\pi,\pi\rrceil_{\varphi^!TM}=0$ and $\Ad_{\aap}(\pi) =\pi$. A Hom-Poisson manifold is a manifold $M$ equipped with a Hom-Poisson tensor $\pi$. We denote a Hom-Poisson manifold by $(M,\varphi,\pi)$.
\end{defi}

In the sequel, we investigate the relation between a Hom-Poisson manifold and a purely Hom-Poisson algebra. Let us recall the definition of a purely Hom-Poisson algebra given in \cite{LGT} first.

\begin{defi}{\rm(\cite{LGT})}
A purely Hom-Poisson algebra is a quadruple $(\frkA,\mu,\{\cdot,\cdot\},\phi)$ consisting of a vector space $\frkA$, bilinear maps $\mu:\frkA\otimes \frkA\rightarrow \frkA$ and $\{\cdot,\cdot\}:\frkA\otimes \frkA\rightarrow \frkA$ and a linear map $\phi: \frkA\rightarrow \frkA$ such that:
\begin{itemize}\item[\rm(i)]
$(\frkA,\mu)$ is a commutative associative algebra,
\item[\rm(ii)]
$(\frkA,\{\cdot,\cdot\},\phi)$ is a Hom-Lie algebra,
\item[\rm(iii)]
$\{x,\mu(y,z)\}=\mu(\phi(y),\{x,z\})+\mu(\{x,y\},\phi(z)),$ for all $x,y,z\in \frkA.$
\end{itemize}
\end{defi}
\begin{rmk}
  The notion of a Hom-Poisson algebra was introduced in \cite{MS2} in the study of formal deformations of a Hom-associative algebra. One difference between a purely Hom-Poisson algebra and a Hom-Poisson algebra is that the former is a commutative associative algebra, whereas the latter is a commutative Hom-associative algebra.
\end{rmk}
\begin{ex}{\rm
  Let $(\g,[\cdot,\cdot]_\g,\phi_\g)$ be a Hom-Lie algebra. Then $\varphi=\phi_\g^*:\g^*\longrightarrow\g^*$  induces an algebra homomorphism $\varphi^*:C^\infty(\g^*)\longrightarrow C^\infty(\g^*)$. In particular, for a linear function $x\in\g$ on $\g^*$, we have $\varphi^*(x)=\phi_\g(x).$ Now define a bilinear map $\{\cdot,\cdot\}:\wedge^2C^\infty(\g^*)\longrightarrow C^\infty(\g^*)$ by
  $$
  \{f,g\}(\xi)=\langle\xi,[\dM f( \xi), \dM g( \xi)]_\g\rangle,\quad \forall ~f,g\in C^\infty(\g^*),~\xi\in\g^*,
  $$
 where $\dM$ is the differential operator of the Hom-Lie algebroid $(\phi_\g^*)^!T\g^*$. Then $(C^\infty(\g^*),\mu,\{\cdot,\cdot\},\varphi^*)$ is a purely Hom-Poisson algebra, where $\mu$ is the usual multiplication on smooth functions. This answers a question about how to associate a purely Hom-Poisson algebra on $S(\g)$ proposed in \cite[Example 2.10]{LGT}.
  }
\end{ex}

Let $M$ be a smooth manifold and   $(\CWM,\{\cdot,\cdot\},\varphi^*)$  a purely Hom-Poisson algebra on $\CWM$. The Hom-Leibniz rule says that $\{f,\cdot\}$ is a $(\varphi^*,\varphi^*)$-derivation on $\CWM$. Thus, there exists a bisection $\pi\in\Gamma(\wedge^2 \varphi^!TM)$ such that
\begin{equation}\label{PHomA}
\{f,g\}=\pi(\dM f,\dM g),\quad f,g\in\CWM.
\end{equation}

By computation, we can obtain the following lemma.
\begin{lem}\label{lem:Hompmor}
With the above notions, the following conditions are equivalent:
\begin{eqnarray}
\varphi^*\{f,g\}&=&\{\varphi^*(f),\varphi^*(g)\};\\
\Ad_{\aap}\pi&=&\pi;\\
\Ad_{\aap}\circ \pi^\sharp&=&\pi^\sharp\circ \Ad^\dag_{\aap},
\end{eqnarray}
where $\Ad^\dag_{\aap}:\Gamma(\varphi^!T^*M)\longrightarrow\Gamma(\varphi^!T^*M)$ is  given by \eqref{eq:auto2} and $\pi^\sharp:\varphi^!T^*M\longrightarrow \varphi^!TM$ is a bundle map defined by
\begin{equation}
\pi^\sharp(\xi)(\eta)=\pi(\xi,\eta),\quad \xi,\eta\in\Gamma(\varphi^!T^*M).
\end{equation}
\end{lem}

\begin{pro}
If $(\CWM,\{\cdot,\cdot\},\varphi^*)$ is a purely Hom-Poisson algebra, then $(M,\varphi,\pi)$ is a Hom-Poisson manifold, where $\pi$ is given by \eqref{PHomA}.

Conversely, If $(M,\varphi,\pi)$ is a Hom-Poisson manifold, then $(\CWM,\{\cdot,\cdot\},\varphi^*)$ is a purely Hom-Poisson algebra,  where $\{\cdot,\cdot\}$ is given by \eqref{PHomA}.
\end{pro}
\pf
By \eqref{eq:HomPFormu}, we have
\begin{eqnarray*}
i_{\llceil\pi,\pi\rrceil_{\varphi^!TM}}(\dM f\wedge\dM g\wedge\dM h)=-2(\varphi^*)^2\Big( \{ (\varphi^*)^{-1} \{ (\varphi^*)^{-1}(f),(\varphi^*)^{-1}(g)\},(\varphi^*)^{-1}(h)\}+c.p.\Big).
\end{eqnarray*}
Since $\varphi^*\{f,g\}=\{\varphi^*(f),\varphi^*(g)\}$, we have
$$i_{\llceil\pi,\pi\rrceil_{\varphi^!TM}}(\dM f\wedge\dM g\wedge\dM h)=-2(\{\{f,g\},\varphi^*(h)\}+\{\{h,f\},\varphi^*(g)\}+\{\{g,h\},\varphi^*(f)\}),$$
which implies that $\llceil\pi,\pi\rrceil_{\varphi^!TM}=0$.
By Lemma \ref{lem:Hompmor}, $(M,\varphi,\pi)$ is a Hom-Poisson manifold.

The converse part can be proved similarly. We omit details.  \qed
\vspace{3mm}

Let $(M,\varphi,\pi)$ be a Hom-Poisson manifold, define a bracket on $\Gamma(\varphi^!T^*M)$ by
\begin{eqnarray}\label{eq:cotangent}
[\xi,\eta]_{\pi^\sharp}=L_{\pi^\sharp(\xi)}\eta-L_{\pi^\sharp(\eta)}\xi-\dM \pi(\xi,\eta),\quad\forall \xi,\eta\in\Gamma(\varphi^!T^*M).
\end{eqnarray}

It is well-known that a Poisson manifold gives rise to a Lie algebroid structure on its cotangent bundle (\cite{cotangentLA}). Similarly, we have

\begin{thm}\label{pro:cotangentHL}
Let $(M,\varphi,\pi)$ be a Hom-Poisson manifold. Then $(\varphi^!T^*M,\varphi,\Ad^\dag_{\varphi^*},[\cdot,\cdot]_{\pi^\sharp},\pi^\sharp)$ is Hom-Lie algebroid, where $[\cdot,\cdot]_{\pi^\sharp}$ is given by \eqref{eq:cotangent}. We call it the {\bf cotangent Hom-Lie algebroid} and denote it by $\varphi_{\pi}^!T^*M$.
\end{thm}
 \pf
For all $f\in\CWM$ and $\xi,\eta\in\Gamma(\varphi^!T^*M)$, since $\Ad_{\aap}\circ \pi^\sharp=\pi^\sharp\circ \Ad^\dag_{\aap}$, we have
\begin{eqnarray}
\nonumber[\xi,f\eta]_{\pi^\sharp}
&=&L_{\pi^\sharp(\xi)}f\eta-L_{\pi^\sharp(f\eta)}\xi-\dM (f\pi(\xi,\eta))\\
\nonumber
&=&\varphi^*(f)[\xi,\eta]_{\pi^\sharp}+\Ad_{\varphi^*}(\pi^\sharp(\xi))(f)\Ad^\dag_{\varphi^*}\eta-\dM f\wedge i_{\pi^\sharp(\eta)}\xi-\dM f\wedge \varphi^*(\pi(\xi,\eta))\\
\nonumber
&=&\varphi^*(f)[\xi,\eta]_{\pi^\sharp}+\Ad_{\varphi^*}(\pi^\sharp(\xi))(f)\Ad^\dag_{\varphi^*}\eta\\
\label{eq:anchor}
&=&\varphi^*(f)[\xi,\eta]_{\pi^\sharp}+\pi^\sharp(\Ad^\dag_{\varphi^*}(\xi))(f)\Ad^\dag_{\varphi^*}\eta.
\end{eqnarray}

For $\xi,\eta,\zeta\in\Gamma(\varphi^!T^*M)$, define
\begin{eqnarray*}
A(\xi,\eta)&=&\pi^\sharp([\xi,\eta]_{\pi^\sharp})\circ \varphi^*-\pi^\sharp(\Ad^\dag_{\aap}\xi)\pi^\sharp(\eta)-\pi^\sharp(\Ad^\dag_{\aap}\eta)\pi^\sharp(\xi),\\
B(\xi,\eta,\zeta)&=&[[\xi,\eta]_{\pi^\sharp},\Ad^\dag_{\aap}(\zeta)]_{\pi^\sharp}+[[\zeta,\xi]_{\pi^\sharp},\Ad^\dag_{\aap}(\eta)]_{\pi^\sharp}+[[\eta,\zeta]_{\pi^\sharp},\Ad^\dag_{\aap}(\xi)]_{\pi^\sharp}.
\end{eqnarray*}
It follows from \eqref{eq:anchor} that $A$ is $\varphi^*$-function linear and $B$ is $(\varphi^*)^2$-function linear.

By \eqref{d0} and \eqref{eq:Lied2}, we have
\begin{eqnarray*}
[\Ad^\dag_{\aap}(\xi),\Ad^\dag_{\aap}(\eta)]_{\pi^\sharp}&=&L_{\pi^\sharp(\Ad^\dag_{\aap}(\xi))}\Ad^\dag_{\aap}(\eta)-L_{\pi^\sharp(\Ad^\dag_{\aap}(\eta))}\Ad^\dag_{\aap}(\xi)-\dM \pi(\Ad^\dag_{\aap}(\xi),\Ad^\dag_{\aap}(\eta))\\
&=&\Ad^\dag_{\aap}(L_{\pi^\sharp(\xi)}\eta-L_{\pi^\sharp(\eta)}\xi)-\dM \varphi^*(\pi(\xi,\eta))\\
&=&\Ad^\dag_{\aap}([\xi,\eta]_{\pi^\sharp}),
\end{eqnarray*}
which implies that
$\Ad^\dag_{\aap}$ is an automorphism.

Since $[\dM f,\dM g]_{\pi^\sharp}=\dM\{f,g\}$, we have
\begin{eqnarray*}
&&[[\dM f,\dM g]_{\pi^\sharp},\Ad^\dag_{\aap}(\dM h)]_{\pi^\sharp}+[[\dM h,\dM f]_{\pi^\sharp},\Ad^\dag_{\aap}(\dM g)]_{\pi^\sharp}+[[\dM g,\dM h]_{\pi^\sharp},\Ad^\dag_{\aap}(\dM f)]_{\pi^\sharp}\\
&=&\{\{f,g\},\varphi^*(h)\}+\{\{h,f\},\varphi^*(g)\}+\{\{g,h\},\varphi^*(f)\}=0,
\end{eqnarray*}
which implies that $B(\dM f,\dM g,\dM h)=0$ and thus $B=0$. Therefore, $(\Gamma(\varphi^!T^*M),[\cdot,\cdot]_{\pi^\sharp},\Ad^\dag_{\varphi^*})$ is a Hom-Lie algebra.

Because of the Hom-Jacobi identity of $\{\cdot,\cdot\}$, we have
\begin{eqnarray*}
&&\pi^\sharp\dM\{f,g\}( \varphi^*(h))
=\{\{f,g\},\varphi^*(h)\}
=\{\varphi^*(f),\{g,h\}\}-\{\varphi^*(g),\{f,h\}\}\\
&=&\varphi^*\circ\pi^\sharp(\dM f)\{(\aap)^{-1}(g),(\aap)^{-1}(h)\}
-\varphi^*\circ\pi^\sharp(\dM g)\{(\aap)^{-1}(f),(\aap)^{-1}(h)\}\\
&=&\varphi^*\circ\pi^\sharp(\dM f)\circ(\aap)^{-1}\circ\pi^\sharp(\dM g)(h)
-\varphi^*\circ\pi^\sharp(\dM g)\circ(\aap)^{-1}\circ\pi^\sharp(\dM f)(h)\\
&=&(\Ad_{\aap}(\pi^\sharp(\dM f))\pi^\sharp(\dM g)-\Ad_{\aap}(\pi^\sharp(\dM g))\pi^\sharp(\dM f))(h)\\
&=&(\pi^\sharp(\Ad^\dag_{\aap}\dM f)\pi^\sharp(\dM g)-\pi^\sharp(\Ad^\dag_{\aap}\dM g)\pi^\sharp(\dM f))(h),
\end{eqnarray*}
which implies that $A(\dM f,\dM g)=0$ and thus $A=0$.

 On the other hand, since $\Ad_{\aap}\circ \pi^\sharp=\pi^\sharp\circ \Ad^\dag_{\aap}$, we have
\begin{eqnarray*}
\pi^\sharp(\Ad^\dag_{\varphi^*}\xi)\circ\varphi^*=\Ad_{\aap}\circ \pi^\sharp(\xi)\circ\varphi^*=\varphi^*\circ\pi^\sharp(\xi).
\end{eqnarray*}
Thus $\pi^\sharp$ is a representation of $(\Gamma(\varphi^!T^*M),[\cdot,\cdot]_{\pi^\sharp},\Ad^\dag_{\varphi^*})$ on $\CWM$ with respect to $\varphi^*$.

Therefore,  $(\varphi^!T^*M,\varphi,\Ad^\dag_{\varphi^*},[\cdot,\cdot]_{\pi^\sharp},\pi^\sharp)$ is Hom-Lie algebroid.
\qed

\section{Hom-Lie bialgebroids}

In this section, we introduce the notion of a Hom-Lie bialgebroid and show that the base manifold of a Hom-Lie bialgebroid is a Hom-Poisson manifold.

\begin{defi}
Let $(\aab,\varphi,\aaa)$ be an invertible Hom-bundle, $(\aab,\varphi,\aaa,[\cdot,\cdot]_{\aab},a_A)$ and $(\aab^*,\varphi,\aad,\\~[\cdot,\cdot]_{\aab^*},a_{A^*})$  two Hom-Lie algebroids in duality. We call $(\aab,\aab^*)$ a Hom-Lie bialgebroid if
\begin{eqnarray}\label{liebi}
\dM_*[ x,y]_\aab=\llceil\dM_*x,\aaa (y)\rrceil_\aab+\llceil\aaa (x),\dM_*y\rrceil_\aab,\quad \forall ~x,y\in\aar,
\end{eqnarray}
where $\dM_*$ is the coboundary operator given by \eqref{D1} for the Hom-Lie algebroid $(\aab^*,\varphi,\aad,[\cdot,\cdot]_{\aab^*},a_{A^*})$.
\end{defi}

\begin{rmk}
When $\varphi=\idd$ and $\phi_A=\idd$, a Hom-Lie algebroid is exactly  a Lie bialgebroid introduced in \cite{Lie bialgebroid}. When the base manifold is a point, a Hom-Lie bialgebroid reduces to a purely Hom-Lie bialgebra given in Definition \ref{defi:Hom-Liebi}.
\end{rmk}

\begin{pro}
Let $(M,\varphi,\pi)$ be a Hom-Poisson manifold. Then $(\varphi^!TM,\varphi_{\pi}^!T^*M)$ is a Hom-Lie bialgebroid, where   $\varphi_{\pi}^!T^*M$ is the cotangent Hom-Lie algebroid associated to Hom-Poisson manifold $(M,\varphi,\pi)$ given in Theorem \ref{pro:cotangentHL}.
\end{pro}
\pf  First consider the differential operator $\dM_*$ associated to the Hom-Lie algebroid  $\varphi_{\pi}^!T^*M$, we have
\begin{eqnarray}
\dM_* X=\llceil \Ad^{-1}_{\aap}(\pi),X\rrceil_{\varphi^*} =\llceil \pi,X\rrceil_{\varphi^*} =\llceil \Ad_{\aap}(\pi),X\rrceil_{\varphi^*} ,\quad X\in\Gamma(\wedge^k A).
\end{eqnarray}
Then by direct calculation, for $x,y\in\Gamma(\varphi^!TM)$, we have
\begin{eqnarray*}
\dM_*[ x,y]_{\varphi^*}
&=&\llceil \Ad_{\aap}(\pi),[ x,y]_{\varphi^*} \rrceil_{\varphi^*} \\
&=&-\llceil \Ad_{\aap}(y),\llceil \pi,x\rrceil_{\varphi^*} \rrceil_{\varphi^*} +\llceil \Ad_{\aap}(x),\llceil \pi,y\rrceil_{\varphi^*} \rrceil_{\varphi^*} \\
&=&\llceil \dM_*x,\Ad_{\aap}(y)\rrceil_{\varphi^*} +\llceil \Ad_{\aap}(x),\dM_*y\rrceil_{\varphi^*},
\end{eqnarray*}
which implies that $(\varphi^!TM,\varphi_{\pi}^!T^*M)$ is a Hom-Lie bialgebroid.\qed

\begin{lem}\label{lll}
Assume that $(\aab,\aab^*)$ is a Hom-Lie bialgebroid. Then we have
\begin{eqnarray}\label{eq:important}
\label{L24}\huaL_{\dM f}x=[x,\dM_*f]_A,\quad\forall~f\in C^{\infty}(M), x\in\aar,
\end{eqnarray}
where $\huaL$ is the Lie derivative of the Hom-Lie algebroid $(\aab^*,\varphi,\aad,[\cdot,\cdot]_{\aab^*},a_{A^*})$
\end{lem}
\pf
For all $f\in C^{\infty}(M)$ and $x,y\in\aar$, by \eqref{liebi}, we have
\begin{eqnarray*}
\dM_*[x,fy]_A
&=&\dM_*\Big{(}\aap(f)[x,y]_A +a_A(\aaa (x))(f)\aaa (y)\Big{)}\\
&=&\dM_*\aap(f)\wedge\aaa([x,y]_A )+(\aap)^2(f)\Big{(}-L_{\aaa (y)}\dM_*x+L_{\aaa (x)}\dM_*y\Big{)}\\
&&+\dM_*a_A(\aaa (x))(f)\wedge\aaa^2(y)+\aap a_A(\aaa (x))(f)\dM_*\aaa (y).
\end{eqnarray*}
On the other hand, we have
\begin{eqnarray*}
\dM_*[x,fy]_A
&=&-L_{\aaa (fy)}\dM_*x+L_{\aaa (x)}\dM_*(fy)\\
&=&-L_{\aaa (fy)}\dM_*x+(L_{\aaa (x)}\dM_* f)\wedge \aaa^2 (y)+\aaa (\dM_* f)\wedge L_{\aaa (x)}(\aaa(y))\\
&&+(\aap)^2(f) L_{\aaa (x)}\dM_*y+a_A(\aaa^2 (x))(\aap(f))\aaa(\dM_*y).
\end{eqnarray*}
By \eqref{eq:LieD4}, we have
 $$L_{\aaa (fy)}\dM_*x=(\aap)^2(f)L_{\aaa (y)}\dM_*x-\aaa^2(y)\wedge i_{\dM f}\dM_*x.$$

It thus follows that
$$(L_{\aaa (x)}\dM_* f)\wedge \aaa^2 (y)=(i_{\dM f}\dM_*x+\dM_*a_A(\aaa (x))(f))\wedge\aaa^2(y),$$
which implies that
$$L_{\aaa (x)}\dM_* f=\huaL_{\dM f}\aaa(x).$$
Thus, \eqref{L24} follows immediately.
 \qed\vspace{3mm}

Denote by $a_A^*:\varphi^!T^*M\rightarrow \aab^*$ and $a_{A^*}^*:\varphi^!T^*M\rightarrow\aab$ the dual map of $a_A$ and $a_{A^*}$ respectively. Define
$\pi^\sharp:=a_A\circ  a_{A^*}^*:\varphi^!T^*M\rightarrow \varphi^!TM.$ We have the following formula.

\begin{cor}\label{lfg}
With the above notions, we have
\begin{eqnarray}\label{76}
[\dM_*f,\dM_*g]_A=\dM_*(\pi^\sharp(\delta f)g).
\end{eqnarray}
where $\delta$ is the differential operator associated to the tangent Hom-Lie algebroid $\varphi^!TM$.
\end{cor}

\begin{thm}\label{thm:Hom-LP}
Let $(\aab,\aab^*)$ be a Hom-Lie bialgebroid. Then $\pi^\sharp:\varphi^!T^*M\rightarrow \varphi^!TM$ defines a Hom-Poisson structure on $M$, and so does $\bar{\pi}^\sharp=a_{A\ast}\circ a_A^\ast:\varphi^!T^*M\rightarrow \varphi^!TM$. Moreover, $\pi^\sharp$ and $\bar{\pi}^\sharp$ are opposite to one another.
\end{thm}
\pf
It is not hard to see that $\pi$ is skew-symmetric. Let $\{f,g\}=\pi^\sharp(\delta f)(g)$ as usual. Then we have
\begin{eqnarray*}
\aap\{f,g\}=\aap\langle\dM f,\dM_*g\rangle=\langle\aad(\dM f),\aaa(\dM_*g)\rangle=\langle\dM\aap(f),\dM_*\aap(g)\rangle=\{\aap(f),\aap(g)\},
\end{eqnarray*}
which implies that $\aap$ is an algebra homomorphism. Applying $a_A$ to both sides of \eqref{76}, we have
$$\pi^\sharp\delta(\{f,g\})=[\pi^\sharp\delta f,\pi^\sharp\delta g]_{\varphi^*},$$
which is equivalent to the Hom-Jacobi identity of the bracket $\{\cdot,\cdot\}$. Thus $\pi$ defines a Hom-Poisson structure on $M$. The rest of the proposition can be obtained directly.
\qed

\begin{pro}\label{equivalent}
If $(\aab,\aab^*)$ is a Hom-Lie bialgebroid, then so is $(\aab^*,\aab)$.
\end{pro}
\pf By direct calculation, we have
\begin{eqnarray*}
&&\langle-\dM[\xi,\eta]_{\aab^*}+\llceil  \aad(\xi),\dM\eta\rrceil  _{\aab^*}-\llceil  \aad(\eta),\dM\xi\rrceil  _{\aab^*},\aaa^2(x)\wedge \aaa^2(y)\rangle\\
&=&\langle-\dM_*{[x,y]_{\aab}}-\llceil  \aaa(y),\dM_*x\rrceil  _{\aab}+\llceil  \aaa(x),\dM_*y\rrceil  _{\aab},(\aad)^2(\xi)\wedge(\aad)^2(\eta)\rangle\\
&&+\huaL_{\dM (i_x\xi)}(i_y\eta)-\huaL_{\dM (i_y\xi)}(i_x\eta)-\huaL_{\dM (i_x\eta)}(i_y\xi)+\huaL_{\dM (i_y\eta)}(i_x\xi)\\
&=&\langle-\dM_*{[x,y]_{\aab}}-\llceil  \aaa(y),\dM_*x\rrceil  _{\aab}+\llceil  \aaa(x),\dM_*y\rrceil  _{\aab},(\aad)^2(\xi)\wedge(\aad)^2(\eta)\rangle\\
&=&0.
\end{eqnarray*}
The second equality follows from $\huaL_{\dM f}(g)=-\huaL_{\dM g}(f)$, which can be obtained from Theorem \ref{thm:Hom-LP}.\qed

\section{Hom-Courant algebroids }

In this section, we introduce the notion of a Hom-Courant algebroid and show that on the double of a Hom-Lie bialgebroid, there is a Hom-Courant algebroid structure. Moreover, we also give the underlying algebraic structure of a Hom-Courant algebroid. First we recall that a {\bf Hom-Leibniz algebra},  which was introduced in \cite{MS2}, is
  a triple $(\frkg,\odot,\phi_\g)$ consisting of a
  vector space $\frkg$, a   bilinear map  $\odot: \frkg\times\g \longrightarrow
  \frkg$ and a linear transformation $\phi_\g:\frkg\lon\frkg$ satisfying $\phi_\g(x\odot y)=\phi_\g(x)\odot\phi_\g(y)$, and the following Hom-Leibniz rule
  \begin{equation}
    \phi_\g(x)\odot(y\odot z)=(x\odot y)\odot\phi_\g(z)+ \phi_\g(y)\odot(x\odot z),\quad\forall
x,y,z\in\frkg.
  \end{equation}

\begin{defi}\label{homcourant}
A Hom-Courant algebroid is an invertible Hom-bundle $(\aae\rightarrow M,\varphi,\aee)$ together with a nondegenerate symmetric bilinear form B on the bundle, a bilinear operation $\odot$ on $\aaee$ and a bundle map $\rho:\aae\rightarrow \varphi^!TM$ such that the following conditions are satisfied:
\begin{itemize}\item[\rm(i)]
$(\aaee,\odot,\aee)$ is a Hom-Leibniz algebra;
\item[\rm(ii)]
$\rho(\aee (e))\circ\varphi^*= \varphi^*\circ \rho(e),\quad \forall e\in\aaee;$
\item[\rm(iii)]
$\rho(e_1\odot e_2)=[\rho(e_1),\rho(e_2)]_{\varphi^*},\quad\forall e_1,e_2\in\aaee;$
\item[\rm(iv)]
$\label{hca2}e\odot e=\frac{1}{2}\ddd B(e,e),~\forall e\in\aaee;$
\item[\rm(v)]
$B(\aee(e_1),\aee(e_2))=\varphi^* B(e_1,e_2),\quad \forall e_1,e_2\in\aaee$;
\item[\rm(vi)]
$\label{hca1}\rho(\aee (e))B(h_1,h_2)=B(e\odot h_1,\aee (h_2))+B(\aee (h_1),e\odot h_2),\quad\forall e,h_1,h_2\in\aaee,$
\end{itemize}
where $\ddd:C^{\infty}(M)\rightarrow\Gamma(\aae)$ is defined by
\begin{equation}
\label{hca22}B(\ddd f,e)=\rho(e)f.
\end{equation}
We denote a Hom-Courant algebroid by $(E,\varphi,\aee,B,\odot,\rho)$.
\end{defi}

\begin{rmk}
  When $\varphi=\idd$ and $\phi_E=\idd$, a Hom-Courant algebroid is exactly a Courant algebroid introduced in  \cite{lwx,Roytenberg4}. When $M$ is a point, a Hom-Courant algebroid reduces to a quadratic Hom-Lie algebra introduced in \cite{CS-bialgebra}.
\end{rmk}

\begin{lem}
Let $(\aae,\varphi,\aee,B,\odot,\rho)$ be a Hom-Courant algebroid. Then for all  $e\in\aaee, f,g\in C^{\infty}(M),$ we have
\begin{eqnarray*}
\label{rho}\rho\circ\ddd&=&0,\\
\label{eq:homlie21}\phi_E \circ D&=&D\circ\varphi^*,\\
\label{dle}\ddd(fg)&=&\ddd(f)\aap(g)+\aap(f)\ddd(g),\\
\label{edf}e\odot\ddd f&=&\ddd B(\ddd f,e),\\
\label{dfe}\ddd f\odot e&=&0.
\end{eqnarray*}
\end{lem}
\pf We only give the proof of the fourth equality, others can be proved similarly.
For $h\in\aaee$, by conditions (iii) and (vi) in Definition \ref{homcourant},
\begin{eqnarray*}
\label{ee1}\rho(\aee (e))(\rho(\aee (h))f)
&=&\rho(\aee (e))(B(\ddd f,\aee (h)))\\
&=&B( e\odot\ddd f,\aee^2(h))+B(\aee(\ddd f),e\odot\aee (h))\\
&=&B( e\odot\ddd f,\aee^2(h))+\rho(e\odot\aee (h))(\aap(f))\\
&=&B( e\odot\ddd f,\aee^2(h))+\rho(\aee (e))(\rho(\aee (h))(f))-\rho(\aee^2(h))(\rho(e)(f)).
\end{eqnarray*}
Hence,
\begin{eqnarray*}
B( e\odot\ddd f,\aee^2(h))=\rho(\aee^2(h))(\rho(e)(f))=B(\ddd\rho(e)f,\aee^2(h))=B(\ddd B(\ddd f,e),\aee^2(h)),
\end{eqnarray*}
which implies that $e\odot\ddd f=\ddd B(\ddd f,e)$.\qed

\begin{lem}
Let $(\aae,\varphi,\aee,B,\odot,\rho)$ be a Hom-Courant algebroid. For all $f\in C^{\infty}(M), e,h\in\aaee,$ we have
\begin{eqnarray}
 \label{xfy}e\odot fh&=&\aap(f) e\odot h+\rho(\aee (e))(f)\aee(h),\\
 \label{xgy}(fe)\odot h&=&\aap(f) e\odot h-\rho(\phi_E(h))(f)\phi_E(e)+\ddd(f)\varphi^*(B(e,h)).
\end{eqnarray}
\end{lem}
\pf By condition (vi) in Definition \ref{homcourant}, we have
\begin{eqnarray*}
\rho(\aee (e)) B( fh_1,h_2)&=& B( e\odot fh_1,\aee (h_2))+ B( \aee(fh_1),e\odot h_2)\\
&=& B( e\odot fh_1,\aee (h_2))+\aap(f) B( \aee (h_1),e\odot h_2).
\end{eqnarray*}
On the other hand, we have
\begin{eqnarray*}
\rho(\aee (e)) B( fh_1,h_2)
&=&\rho(\aee (e))(f B( h_1,h_2))\\
&=&\rho(\aee (e))(f)\aap( B( h_1,h_2))+\aap(f)\rho(\aee (e))( B( h_1,h_2))\\
&=&\rho(\aee (e))(f)\aap( B( h_1,h_2))
+\aap(f) B( e\odot h_1,\aee (h_2))\\&&+\aap(f) B( \aee (h_1),e\odot h_2).
\end{eqnarray*}
Thus, we have
\begin{eqnarray*}
 B( e\odot fh_1,\aee (h_2))&=&\aap(f) B( e\odot h_1,\aee (h_2))+\rho(\aee (e))(f)\aap( B( h_1,h_2))\\
&=&\aap(f) B( e\odot h_1,\aee (h_2))+\rho(\aee (e))(f) B( \aee(h_1),\aee(h_2))\\
&=& B(\aap(f) e\odot h_1+\rho(\aee (e))(f)\aee(h_1),\aee(h_2)),
\end{eqnarray*}
which implies that \eqref{xfy} holds.\qed\vspace{3mm}

Suppose that $(A,\varphi,\phi_A,[\cdot,\cdot]_A,a_A)$ and $(A^*,\varphi,\aad,[\cdot,\cdot]_{A^*},a_{A^*})$ are Hom-Lie algebroids over the base manifold $M$. Let $E$ denote their vector bundle direct sum: $E=A\oplus A^*$. On $\Gamma(E)$, there is a natural nondegenerate symmetric bilinear form given by
\begin{equation}\label{eq:pair}
   B( x+\xi,y+\eta)=\xi(y)+\eta(x),\quad \forall ~x, y\in\Gamma(A),\xi, \eta\in\Gamma(A^*).
\end{equation}
On $\Gamma(E)$,  we introduce an operation $\odot:\Gamma(E)\times\Gamma(E)\longrightarrow\Gamma(E)$ by
\begin{eqnarray}\label{A2}
(x+\xi)\odot(y+\eta)
=\big{(}[x,y]_\aab+\huaL_\xi y-i_\eta\dM_*\aaa^{-1}(x)\big{)}
+\big{(}[\xi,\eta]_{\aab^*}+L_x\eta-i_y\dM\aadd(\xi)\big{)}.
\end{eqnarray}
Define $\aee:\Gamma(E)\longrightarrow \Gamma(E)$ by $\aee=\aaa\oplus\aad$. That is
\begin{equation}
\label{eq:HCAauto} \aee(x+\xi)=\aaa(x)+\aad(\xi).
\end{equation}
Finally, we let $\rho:\aae\rightarrow \varphi^!TM$ be the bundle map defined by $\rho=a_A+a_{A^*}$. That is
\begin{equation}
\label{eq:HCAanchor} \rho(x+\xi)=a_A(x)+a_{A^*}(\xi).
\end{equation}

\begin{thm}
Let $(\aab,\aab^*)$ be a Hom-Lie bialgebroid. Then
$(E=\aab\oplus\aab^*,\varphi,\aee, B,\odot,\rho)$ is a Hom-Courant algebroid, where $ B$ is given by \eqref{eq:pair}, $\odot$ is given by \eqref{A2}, $\aee$ is given by \eqref{eq:HCAauto} and $\rho$ is given by \eqref{eq:HCAanchor}.
\end{thm}

\pf
Since $a_A$ and $a_{A^*}$ are representations of Hom-Lie algebroids $A$ and $A^*$ respectively, it is easy to see that condition (ii) in  Definition \ref{homcourant} holds.

By \eqref{eq:important}, for $f\in\CWM$, we have
\begin{eqnarray*}
\rho(x\odot\eta)(f)&=&-a_A(i_\eta\dM_*\aaa^{-1}(x))(f)+a_{A^*}(L_x\eta)(f)
=-(\dM_*x)(\aad(\eta),\dM f)+(L_x\eta)(\dM_*f)\\
&=&-a_{A^*}(\aad(\eta))(x(\dM(\aap)^{-1}(f)))+a_A(\dM f)(x(\eta))+\aaa (x)([\eta,\dM(\aap)^{-1}(f)]_{\aab^*})\\
&&+a_A(\aaa (x))(\eta(\dM_*(\aap)^{-1}(f)))-(\aad(\eta))([x,\dM_*(\aap)^{-1}(f)]_\aab)\\
&=&-a_{A^*}(\aad(\eta))\Big{(}a_A(x)((\aap)^{-1}(f))\Big{)}+(\huaL_{\dM(\aap)^{-1}(f)}x)(\aad(\eta))\\
&&+a_A(\aaa (x))\Big{(}a_{A^*}(\eta)((\aap)^{-1}(f))\Big{)}-(\aad(\eta))([x,\dM_*(\aap)^{-1}(f)]_\aab)\\
&=&-a_{A^*}(\aad(\eta))\Big{(}a_A(x)((\aap)^{-1}(f))\Big{)}+a_A(\aaa (x))\Big{(}a_{A^*}(\eta)((\aap)^{-1}(f))\Big{)}\\
&=&[\rho(x),\rho(\eta)]_{\varphi^*}(f).
\end{eqnarray*}
Also, for $x,y\in\Gamma(A),\xi,\eta\in\Gamma(A^*)$, we have
$$a_A([x,y]_A)=[a_A(x),a_A(y)]_{\varphi^*},\quad a_{A^*}([\xi,\eta]_{A^*})=[a_{A^*}(\xi),a_{A^*}(\eta)]_{\varphi^*}.$$
Thus, for all $e_1,e_2\in\aaee$, we have
$\rho(e_1\odot e_2)=[\rho(e_1),\rho(e_2)]_{\varphi^*}$, i.e. condition (iii) in  Definition \ref{homcourant} holds.

By direct calculation, conditions (iv), (v) and (vi) in  Definition \ref{homcourant} follow immediately. We omit the detail.

At last, we prove that $(\aaee,\odot,\aee)$ is a Hom-Leibniz algebra.
By \eqref{d0}, \eqref{L20} and \eqref{eq:Lied2}, we can obtain that
$$\aee(e_1\odot e_2)=\aee(e_1)\odot \aee(e_2),\quad \forall e_1,e_2\in\Gamma(E).$$
In the following, we prove that the Hom-Leibniz identity holds, i.e.
\begin{equation}\label{eq:hom-leibniz}
\aee(e_1)\odot (e_2\odot e_3)=(e_1\odot e_2)\odot\aee(e_3)+\aee(e_2)\odot (e_1\odot e_3).
\end{equation}
By \eqref{L20}, \eqref{eq:LieD1} and \eqref{eq:Lied2}, for $x,y\in\Gamma(A),\xi,\eta\in\Gamma(A^*)$, we have
  \begin{eqnarray*}
 && B(\aaa (x)\odot( y\odot\xi)-(x\odot y)\odot\aee(\xi)-\aee(y)\odot (x\odot \xi),\eta)\\
&=& B( i_{\aad(\xi)}\dM_*\aaa^{-1}([x,y]_{\aab})+i_{L_x\xi}\dM_*y-i_{L_y\delta}\dM_*x-[\aaa (x),i_\xi\dM_*\aaa^{-1}(y)]_{\aab}+[\aaa (y),i_\xi\dM_*\aaa^{-1}(x)]_{\aab},\eta)\nonumber\\
&=&(\dM_*[x,y]_{\aab})((\aad)^2(\xi),\eta)+(\dM_*\aaa (y))(L_{\aaa (x)}\aad(\xi),\eta)-(\dM_*\aaa (x))(L_{\aaa (y)}\aad(\xi),\eta)\\
&&+(\dM_*\aaa (y))((\aad)^2(\xi),L_{\aaa (x)}\aadd(\eta))-a_A(\aaa^2(x))((\dM_*y)(\aad(\xi),\aadd(\eta)))\\
&&-(\dM_*\aaa (x))((\aad)^2(\xi),L_{\aaa (y)}\aadd(\eta))+a_A(\aaa^2(y))((\dM_*x)(\aad(\xi),\aadd(\eta)))\\
&=&(\dM_*[x,y]_{\aab})((\aad)^2(\xi),\eta)-((\aad)^2(\xi)\wedge\eta)(\llceil\aaa (x),\dM_*y\rrceil_{\aab})+((\aad)^2(\xi)\wedge\eta)(\llceil\aaa (y),\dM_*x\rrceil_{\aab})\\
&=& B(\dM_*[x,y]_{\aab}-\llceil\aaa (x),\dM_*y\rrceil_{\aab}+\llceil\aaa (y),\dM_*x\rrceil_{\aab},(\aad)^2(\xi)\wedge\eta)\\
&=&0,
 \end{eqnarray*}
 which implies that \eqref{eq:hom-leibniz} holds for $x,y\in\Gamma(A),\xi\in\Gamma(A^*)$.

 By Proposition \ref{equivalent},  \eqref{eq:hom-leibniz} holds for $x\in\Gamma(A),\xi,\eta\in\Gamma(A^*)$. Therefore, $(E=\aab\oplus\aab^*,\varphi, B,\odot,\aee,\rho)$ is a Hom-Courant algebroid.\qed
\begin{rmk}
  When $\varphi=\idd$ and $\phi_A=\idd$, a Hom-Lie bialgebroid is a Lie bialgebroid. The Hom-Courant algebroid given in the above theorem is exactly the Courant algebroid on the double of a Lie bialgebroid (\cite{lwx}). When $M$ is a point, a Hom-Lie bialgebroid reduces to a purely Hom-Lie bialgebra (Definition \ref{defi:Hom-Liebi}). The   Hom-Courant algebroid given in the above theorem is exactly the quadratic Hom-Lie algebra given in Theorem \ref{thm:doublealg}.
\end{rmk}

Associated to any manifold $M$ and a diffeomorphism $\varphi:M\longrightarrow M$, there is a standard Hom-Courant algebroid.
\begin{ex}{\rm
  $(\varphi^!TM\oplus \varphi^!T^*M,\varphi,\Ad_{\varphi^*}\oplus \Ad_{\varphi^*}^\dag,B,\odot,\rho=pr_{\varphi^!TM})$ is a Hom-Courant algebroid, where the symmetric nondegenerate bilinear form $B$ is given by \eqref{eq:pair} for all $x,y\in\Gamma(\varphi^!TM), \xi,\eta\in\Gamma(\varphi^!T^*M)$ and the bilinear map  $\odot$ is given by
  $$
  (x+\xi)\odot(y+\eta)
= [x,y]_\aab+L_x\eta-i_y\dM\aadd(\xi).
  $$
  }
\end{ex}

Finally, we study the algebraic structure underlying of a Hom-Courant algebroid. Similar as the fact that a Courant algebroid gives rise to a Lie $2$-algebra, a Hom-Courant algebroid  gives rise to a Hom-Lie $2$-algebra.
\begin{defi}{\rm(\cite{SC-2-algebra})}\label{defi:2hl}
  A Hom-Lie $2$-algebra consists of the following data:
\begin{itemize}
\item[$\bullet$] a complex of vector spaces $\huaV:V_1\stackrel{l_1}{\longrightarrow}V_0,$

\item[$\bullet$] bilinear maps $l_2,:V_i\times V_j\longrightarrow
V_{i+j}$,

\item[$\bullet$] a skew-symmetric trilinear map $l_3:V_0\times V_0\times V_0\longrightarrow
V_1$,

\item[$\bullet$] two linear transformations $\phi_0\in\gl(V_0)$ and $\phi_1\in\gl(V_1)$ satisfying $$\phi_0\circ l_1=l_1\circ\phi_1,~\phi_0\circ l_2=l_2\circ(\phi_0\times \phi_0), ~\phi_1\circ l_2=l_2\circ(\phi_0\times \phi_1), ~l_3\circ\phi_0=\phi_1\circ l_3$$
   \end{itemize}
   such that for any $w,x,y,z\in V_0$ and $m,n\in V_1$, the following equalities are satisfied:
\begin{itemize}
\item[$\rm(a)$] $l_2(x,y)=-l_2(y,x),$  $l_2(x,m)=-l_2(m,x),$
\item[$\rm(b)$]$l_1 l_2(x,m)=l_2(x,l_1( m)),$  $l_2(l_1( m),n)=l_2(m,l_1( n)),$

\item[$\rm(c_1)$]  $l_1l_3(x,y,z)=l_2(\phi_0(x),l_2(y,z))+l_2(\phi_0(y),l_2(z,x))+l_2(\phi_0(z),l_2(x,y)),$
\item[$\rm(c_2)$] $  l_3(x,y,l_1( m))=l_2(\phi_0(x),l_2(y,m))+l_2(\phi_0(y),l_2(m,x))+l_2(\phi_1(m),l_2(x,y)),$

\item[$\rm(d)$] \begin{eqnarray*}
&&l_3(l_2(w,x),\phi_0(y),\phi_0(z))+l_2(l_3(w,x,z),\phi^2_0(y))\\
&&+
l_3(\phi_0(w),l_2(x,z),\phi_0(y))+l_3(l_2(w,z),\phi_0(x),\phi_0(y)) \\
&=&l_2(l_3(w,x,y),\phi^2_0(z))+l_3(l_2(w,y),\phi_0(x),\phi_0(z))+l_3(\phi_0(w),l_2(x,y),\phi_0(z))\\
&&+l_2(\phi^2_0(w),l_3(x,y,z))+l_2(l_3(w,y,z),\phi^2_0(x))+l_3(\phi_0(w),l_2(y,z),\phi_0(x)).
\end{eqnarray*}
   \end{itemize}
\end{defi}

We  denote a Hom-Lie 2-algebra by
$(\huaV,l_1,l_2,l_3,\phi_0,\phi_1)$.\vspace{3mm}

Let $(E,\varphi, B( \cdot,\cdot),\odot,\aee,\rho)$ be a Hom-Courant algebroid. We introduce a new bracket on $\Gamma(E)$,
$$
\llbracket e_1,e_2\rrbracket=\half(e_1\odot e_2-e_2\odot e_1),\quad e_1,e_2\in\Gamma(E),
$$
which is the skew-symmetrization of $\odot$.

Consider the graded vector space $\huaV=\huaV_0\oplus \huaV_1$, where $V_0=\Gamma(E)$ and $V_1=\CWM$.

\begin{thm}
A Hom-Courant algebroid $(E,\varphi,\aee, B,\odot,\rho)$ gives rise to a Hom-Lie $2$-algebra $(\huaV,l_1,l_2,l_3,\phi_0,\phi_1)$, where $\phi_0,~\phi_1$ and $l_1,~l_2,~l_3$ are given by the following formulas:
\begin{eqnarray*}
\phi_0(e)&=&\phi_E(e),\\
\phi_1(f)&=&\varphi^*(f),\\
l_1(f)&=&\huaD f,\\
l_2(e_1\wedge e_2)&=&\llbracket e_1,e_2\rrbracket,\\
l_2(e\wedge f)&=&\frac{1}{2} B( e,\huaD f),\\
l_3(e_1\wedge e_2\wedge e_3)&=&-T(e_1,e_2,e_3),
\end{eqnarray*}
for all $f\in\CWM,e\in\Gamma(E),e_1,e_2,e_3\in\Gamma(E)$, where $T:\Gamma(E)\times \Gamma(E) \times \Gamma(E)\longrightarrow \CWM$ is given by
$$
T(e_1,e_2,e_3)=\frac{1}{6} B(\llbracket e_1,e_2\rrbracket,\phi_E(e_3))+c.p.,\quad \forall e_1,e_2,e_3\in\Gamma(E).
$$
\end{thm}
\pf The proof is parallel to the proof of \cite[Theorem 4.3]{RW}. We omit details. \qed


\begin{thebibliography}{999}


\bibitem{AEM} F. Ammar, Z. Ejbehi and A. Makhlouf, Cohomology and deformations
of Hom-algebras. \emph{J. Lie Theory} 21 (2011), no. 4, 813-836.




\bibitem{BM}
S. Benayadi and A. Makhlouf, Hom-Lie algebras with symmetric
invariant nondegenerate bilinear forms. \emph{J. Geom. Phys.} 76 (2014), 38-60.


 \bibitem{BEM-quantization}
 M. Bordemann, O. Elchinger and A. Makhlouf,   Twisting Poisson algebras, coPoisson algebras and quantization.
\emph{Trav. Math. } 20 (2012), 83-119.


\bibitem{CS-bialgebra}
L. Cai and Y. Sheng, Purely Hom-Lie bialgebras. arXiv:1605.00722.

\bibitem{cotangentLA}
B. Fuchssteiner, The Lie algebra structure of degenerate Hamiltonian and bi-Hamiltonian systems. \emph{Progr. Theoret. Phys.} 68 (1982), no. 4, 1082-1104.

\bibitem{HLS}
J. Hartwig, D. Larsson and S. Silvestrov, Deformations of Lie
algebras using $\sigma$-derivations. \emph{J. Algebra} 295 (2006), no. 2,
314-361.

\bibitem{Hu}
N. Hu, $q$-Witt algebras, $q$-Lie algebras,
$q$-holomorph structure and representations. \emph{ Algebra Colloq.}
6 (1999), no. 1, 51-70.


\bibitem{Ka}
C. Kassel,
Cyclic homology of differential operators, the Virasoro algebra and a
$q$-analogue, \emph{Comm. Math. Phys.}
146 (1992), no. 2, 343-351.

\bibitem{LGMT}
C. Laurent-Gengoux, A. Makhlouf and J. Teles, Universal algebra of a Hom-Lie algebra and group-like elements.  arXiv:1505.02439.



\bibitem{LGT}
C. Laurent-Gengoux and J. Teles, Hom-Lie algebroids. \emph{J. Geom.  Phys.} 68 (2013), 69-75.


\bibitem{LD1}
D. Larsson and S. Silvestrov,  Quasi-Hom-Lie algebras, central
extensions and 2-cocycle-like identities. \emph{ J. Algebra} 288
(2005), no. 2, 321-344.

\bibitem{LD2} D. Larsson and S. Silvestrov, Graded quasi-Lie algebras. \emph{Czechoslovak J. Phys.} 55 (2005), 11, 1473-1478.


\bibitem{lwx}
Z. Liu, A. Weinstein and P. Xu.
\newblock Manin triples for {L}ie bialgebroids.
\newblock {\em J. Diff. Geom.} 45 (1997), no. 3, 547-574.

\bibitem{General theory of Lie groupoid and Lie algebroid}
K. C. H. Mackenzie, General theory of Lie groupoids and Lie
algebroids.  \emph{Lecture Note Series, $213$. London Mathematical
Society.}, Cambridge University Press, Cambridge, 2005.

\bibitem{Lie bialgebroid} K. C. H. Mackenzie and P. Xu, Lie bialgebroids and Poisson groupoids. \emph{Duke Math. J.} 73 (1994), no. 2, 415-452.

\bibitem{MS2} A. Makhlouf and S. Silvestrov, Hom-algebra
structures. \emph{J. Gen. Lie Theory Appl.}  2 (2008), no. 2,
51-64.

\bibitem{MS1}
 A. Makhlouf and S. Silvestrov, Notes on formal feformations of Hom-associative and Hom-Lie
 algebras. \emph{Forum Math.} 22 (2010), no. 4, 715-739.



\bibitem{Roytenberg4}
D. Roytenberg, \emph{Courant algebroids, derived brackets and even
symplectic supermanifolds}, PhD thesis, UC Berkeley, 1999,
arXiv:math.DG/9910078.

\bibitem{RW}
D. Roytenberg and A. Weinstein, Courant algebroids and strongly homotopy Lie
algebras, \emph{Lett. Math. Phys.} 46 (1998), no. 1, 81-93.



\bibitem{SH-rep}
Y. Sheng, Representations of Hom-Lie algebras. \emph{Algebr. Represent. Theory} 15 (2012), no. 6, 1081-1098.

\bibitem{SC-2-algebra}
Y. Sheng and D. Chen, Hom-Lie 2-algebras. \emph{J. Algebra} 376 (2013) 174-195.

\bibitem{SX}
Y. Sheng and Z. Xiong, On Hom-Lie algebras. \emph{Linear Multilinear Algebra}   63 (2015), no. 12, 2379-2395.

\end{thebibliography}
\end{document}